\documentclass[preprint]{elsarticle}

\usepackage{template}
\usepackage{framed,multirow}

\usepackage{amssymb}
\usepackage{latexsym}
\usepackage{slashbox}

\usepackage{url}
\usepackage{ulem}	
\usepackage{color}
\usepackage{lineno}
\usepackage{amsmath,mathrsfs}
\usepackage{amsthm}
\usepackage[caption=false]{subfig}
\usepackage{dsfont}
\usepackage{float}
\usepackage{gensymb}
\usepackage{natbib}
\usepackage{graphicx,titlesec}
\usepackage{epstopdf}
\usepackage{geometry}
\usepackage{booktabs}
\usepackage{makecell}

\usepackage{varioref}
\usepackage{hyperref}
\usepackage{cleveref}

\hypersetup{
	colorlinks=true,
	linkcolor={red!50!black},
	citecolor={blue!50!black},
	urlcolor={blue!80!black}
}

\newcommand{\abs}[1]{\left\lvert#1\right\rvert}

\newtheorem{remark}{Remark}

\crefformat{equation}{(#2#1#3)}
\crefmultiformat{equation}{(#2#1#3)}{ and~(#2#1#3)}{, (#2#1#3)}{ and~(#2#1#3)}

\crefformat{figure}{Fig.~#2#1#3}
\crefmultiformat{figure}{Figures~(#2#1#3)}{ and~(#2#1#3)}{, (#2#1#3)}{ and~(#2#1#3)}
\crefformat{table}{Table~#2#1#3}

\def\bx{\mbox{\boldmath $x$}}
\def\by{\mbox{\boldmath $y$}}

\graphicspath{{figures/}}

\begin{document}

\verso{Given-name Surname \textit{etal}}

\begin{frontmatter}

\title{A Gauss-Seidel projection method with the minimal number of updates for stray field in micromagnetic simulations}

\author[1]{Panchi Li}
\ead{LiPanchi1994@163.com}
\author[1]{Zetao Ma}
\ead{770120068@qq.com}
\author[1,2]{Rui Du\corref{cor1}}
\ead{durui@suda.edu.cn}
\author[1,2]{Jingrun Chen\corref{cor1}}
\ead{jingrunchen@suda.edu.cn}
\cortext[cor1]{Corresponding authors}

\address[1]{School of Mathematical Sciences, Soochow University, Suzhou, 215006, China.}
\address[2]{Mathematical Center for Interdisciplinary Research, Soochow University, Suzhou, 215006, China.}

\begin{abstract}
Magnetization dynamics in magnetic materials is often modeled by the Landau-Lifshitz equation, which is solved numerically in general. In micromagnetic simulations, the computational cost relies heavily on the time-marching scheme and the evaluation of stray field. Explicit marching schemes are efficient but suffer from severe stability constraints, while nonlinear systems of equations have to be solved in implicit schemes though they are unconditionally stable. A better compromise between stability and efficiency is the semi-implicit scheme, such as the Gauss-Seidel projection method (GSPM) and the second-order backward differentiation formula scheme (BDF2). At each marching step, GSPM solves several linear systems of equations with constant coefficients and updates the stray field several times, while BDF2 updates the stray field only once but solves a larger linear system of equations with variable coefficients and a nonsymmetric structure. In this work, we propose a new method, dubbed as GSPM-BDF2, by combing the advantages of both GSPM and BDF2. Like GSPM, this method is first-order accurate in time and second-order accurate in space, and is unconditionally stable with respect to the damping parameter. However, GSPM-BDF2 updates the stray field only once per time step, leading to an efficiency improvement of about $60\%$ than the state-of-the-art GSPM for micromagnetic simulations. For Standard Problem \#4 and \#5 from National Institute of Standards and Technology, GSPM-BDF2 reduces the computational time over the popular software OOMMF by $82\%$ and $96\%$, respectively. Thus, the proposed method provides a more efficient choice for micromagnetic simulations.
\end{abstract}


\end{frontmatter}




\section{Introduction}
Due to their intrinsic magnetic properties, ferromagnets have been ideal materials for magnetic recording devices in the past several decades \cite{ZuticFabianDasSarma:2004}. The basic quantity of interest in a ferromagnet is the magnetization, whose dynamics is modeled by the Landau-Lifshitz~(LL) equation~\cite{LandauLifshitz1935,Gilbert1955} phenomenologically. The model and its generalizations in the presence of external controls, such as spin current and temperature gradient, have been successfully used to interpret many interesting experimental observations. Numerically, micromagnetic simulation becomes increasingly important to study the dynamics of magnetization, in addition to experiment and theory. In the LL equation, magnetization dynamics is driven by the gyromagnetic term and the damping term, which are both nonlinear with respect to magnetization. Moreover, the length of magnetization does not change during its dynamic evolution. These pose challenges in designing efficient and simple numerical methods. Besides, the calculation of the stray field in micromagnetic simulations is time-consuming since it involves a problem defined over the entire space instead of the ferromagnetic body~\cite{strayfield2013review}.

Explicit schemes such as Runge-Kutta method ~\cite{FourRK2008} are favored in the early days, and the fifth-order Cash-Karp Runge-Kutta method is still used in OOMMF~\cite{fifthRungeKuta1990OOMMF,Najafi2009std}. These schemes are simple and efficient in the sense that no linear/nonlinear systems of equations need to be solved at each marching step, but the temporal stepsize is rather small due to the strong stability restriction of explicit schemes. Implicit methods, such as~\cite{Yamada2004Implicit,bartels2006convergence,CKwithoutdamping2010,implicit2012}, are proved to be unconditionally stable, and can preserve the length of magnetization automatically. However, a nonlinear system of equations with variable coefficients and a nonsymmetric structure has to be solved at each step.

A nice compromise between efficiency and stability is the semi-implicit scheme. One notable example is the Gauss-Seidel projection method~(GSPM)~\cite{NumGSPM2001,ImproveGSPM2003}. At each step, GSPM only solves linear systems of equations with constant coefficients seven times and updates the stray field four times. GSPM is tested to be unconditionally stable, and is first-order accurate in time and second-order accurate in space. Recently, it has been found that the number of linear systems to be solved and the number of updates for the stray field at each step can be reduced without sacrificing the stability and accuracy~\cite{panchi2020GSPM, panchi2020IEEE}. This leads to a reduction in computational cost of about $30\%$ in micromagnetic simulations. Other semi-implicit schemes, such as~\cite{seim2005,Changjian2019semiimplicit,jingrun2019analysis}, are second-order accurate in time. In~\cite{Changjian2019semiimplicit}, the semi-implicit schemes are constructed using the backward differentiation formula (BDF) and one-sided interpolation, and the second-order BDF scheme (BDF2) is proved to converge with second-order accuracy in both space and time~\cite{jingrun2019analysis}. Note that a projection step is needed in semi-implicit schemes to preserve the length of magnetization at each step.

To compare GSPM and BDF2 in details, let us consider a discrete problem with $n$ degrees of freedom (dofs) and the dofs of magnetization is thus $3n$. At each step, GSPM~(referred to Scheme A in \cite{panchi2020GSPM} throughout the paper) solves linear systems of equations with constant coefficients and $n$ dofs  five times and updates the stray field with $3n$ dofs three times. BDF2 solves linear systems of equations once with variable coefficients, a nonsymmetric structure, and $3n$ dofs and updates the stray field once with $3n$ dofs. Both schemes are second-order accurate in space. GSPM is first-order accurate in time while BDF2 is second-order accurate in time. A natural question arises here: \textit{can one design a more efficient method which combines the strengths of both GSPM and BDF2}?

In the work, we provide an affirmative answer to the aforementioned question. The new method is dubbed as GSPM-BDF2, which solves linear systems of equations with constant coefficients and $n$ dofs five times and updates the stray field with $3n$ dofs only once. The method is tested to be unconditionally stable with respect to the damping parameter, first-order accurate in time, and second-order accurate in space. In micromagnetic simulations, GSPM-BDF2 reduces the number of evaluations of the stray field from $3$ to $1$, yielding about $60\%$ reduction of computational time on top of GSPM~\cite{panchi2020GSPM}. For Standard Problem \#4 and \#5 from National Institute of Standards and Technology (NIST)~\cite{NIST}, GSPM-BDF2 reduces the computational time over OOMMF~\cite{OOMMF} by $82\%$ and $96\%$, respectively.

The paper is organized as follows. The LL equation is introduced in \Cref{sec:model}. The GSPM and the proposed method: GSPM-BDF2 are given in \Cref{sec:numericalmethods}. In \Cref{sec:numericalsimulation}, the accuracy and stability of GSPM-BDF2 with respect to the damping parameter are checked in both 1D and 3D, and two benchmark problems from NIST are simulated. Conclusions are drawn in \Cref{sec:conclusion}.

\section{Landau-Lifshitz equation}
\label{sec:model}
The dynamics of magnetization $\mathbf{M} = (M_1, M_2, M_3)^T$ in ferromagnetic materials is modeled by the phenomenological Landau-Lifshitz equation~\cite{LandauLifshitz1935,Gilbert1955}
\begin{equation}
\label{equ:LL_equation}
\mathbf{M}_t = -\gamma\mathbf{M}\times\mathbf{H} -
\frac {\gamma\alpha}{M_s}\mathbf{M}\times(\mathbf{M}\times\mathbf{H}),
\end{equation}
where $\gamma$ is the gyromagnetic ratio, $\alpha$ is the dimensionless damping parameter and $|\mathbf{M}| = M_s$ in the point-wise sense with $M_s$ the saturation magnetization. The local effective field $\mathbf{H} = -\frac{\delta F}{\delta\mathbf{M}}$ is obtained from the Landau-Lifshitz energy functional
\begin{equation}
\label{equ:LL-Energy}
F[\mathbf{M}] = \frac 1{2}\int_\Omega\left\{ \frac{A}{M_s^2}|\nabla\mathbf{M}|^2 + \Phi\left(\frac{\mathbf{M}}{M_s}\right) -2\mu_0\mathbf{H}_e\cdot\mathbf{M} \right\}\mathrm{d}\bx + \frac{\mu_0}{2}\int_{\mathbb{R}^3}|\nabla U|^2 \mathrm{d}\bx
\end{equation}
with $A$ the exchange constant and $\mu_0$ the permeability of vacuum. $\mathbf{H}_e$ is the external magnetic field and $\Omega$ is the volume occupied by the material. For a uniaxial material with $x$-axis the easy direction, the anisotropy energy can be described by $\Phi\left(\frac{\mathbf{M}}{M_s}\right) = \frac{K_u}{M_s^2}(M_2^2+M_3^2)$ with $K_u$ the anisotropy constant. The last term in \eqref{equ:LL-Energy} is the self-interacting energy induced by the magnetization distribution inside the material with
\begin{equation}
U(\bx) = \int_{\Omega}\nabla N(\bx-\by)\cdot\mathbf{M}(\by)\mathrm{d}\by,
\label{equ:scalarU}
\end{equation}
where $N(\bx-\by) = -\frac 1{4\pi}\frac 1{|\bx-\by|}$ is the Newtonian potential. Fast Fourier Transform (FFT) is employed for the evaluation of the stray field $\boldsymbol{\mathbf{H}}_s = -\nabla U$ in regular-shaped materials~\cite{CJreview2007,thesis2008Yanglei}.

For convenience, we nondimensionalize \eqref{equ:LL_equation} by rescaling variables $t \rightarrow(\mu_0\gamma M_s)^{-1}t$ and  $x\rightarrow Lx$ with $L$ the diameter of $\Omega$. With $\mathbf{m}=\mathbf{M}/M_s$, $\mathbf{h}_e = \mathbf{H}_e/M_s$, and $\mathbf{h}_s = \mathbf{H}_s/M_s$, the dimensionless LL equation reads as
\begin{equation}
\label{equ:LLdimensionless}
\mathbf{m}_t = -\mathbf{m}\times(\epsilon\Delta\mathbf{m} + \mathbf{f}(\mathbf{m})) - \alpha\mathbf{m}\times(\mathbf{m}\times(\epsilon\Delta\mathbf{m} + \mathbf{f}(\mathbf{m}))),
\end{equation}
where
\begin{equation}
\mathbf{f}(\mathbf{m}) = -Q(m_2\mathbf{e}_2 + m_3\mathbf{e}_3) + \mathbf{h}_e + \mathbf{h}_s.
\label{equ:effectivefield_withoutSTT}
\end{equation}
Here dimensionless parameters $\epsilon=A/(\mu_0M_s^2L^2)$, $Q=K_u/(\mu_0M_s^2)$, $\mathbf{e}_2=(0,1,0)^T$, and $\mathbf{e}_3=(0,0,1)^T$.
Neumann boundary condition is used
\begin{equation}
\frac{\partial\mathbf{m}}{\partial\boldsymbol{\nu}}\Big\vert_{\partial\Omega} = 0,
\label{equ:Neumannboundary}
\end{equation}
where $\boldsymbol{\nu}$ is the outward unit normal vector on $\partial\Omega$.

An equivalent form of \eqref{equ:LL_equation} is the so-called Landau-Lifshitz-Gilbert (LLG) equation
\begin{equation}
\label{equ:LLGequ}
\mathbf{M}_t = -\gamma\mathbf{M}\times\mathbf{H} + \frac{\alpha}{M_s}\mathbf{M}\times\mathbf{M}_t.
\end{equation}
In the presence of spin transfer torque~(STT), the generalized LLG equation reads as~\cite{ZhangLi2004STT}
\begin{equation}
\mathbf{M}_t = -\gamma\mathbf{M}\times\boldsymbol{\mathbf{H}} + \frac{\alpha}{M_s}\mathbf{M}\times\mathbf{M}_t - \frac{b}{M_s^2}\mathbf{M}\times(\mathbf{M}\times(\mathbf{j}\cdot\nabla)\mathbf{M}) - \frac{b\xi}{M_s}\mathbf{M}\times(\mathbf{j}\cdot\nabla)\mathbf{M},
\label{equ:LLG_STT}
\end{equation}
where $\mathbf{j}$ is the spin polarization current density with magnitude $J$, $b = P\mu_B\left(eM_s(1+\xi^2)\right)^{-1}$ with $e$ the electron charge, $\mu_B$ the Bohr magneton, and $P$ the spin current polarization. After rescaling $t\rightarrow(1+\alpha^2)(\mu_0\gamma M_s)^{-1}t$ and $x\rightarrow Lx$, the LLG equation can be rewritten into \eqref{equ:LLdimensionless} with the dimensionless local field $\mathbf{f}(\mathbf{m})$ of the following form
\[
\mathbf{f}(\mathbf{m}) = -Q(m_2\mathbf{e}_2 + m_3\mathbf{e}_3) + \mathbf{h}_e + \mathbf{h}_s + \frac{b}{\gamma M_sL}\mathbf{m}\times(\mathbf{j}\cdot\nabla)\mathbf{m} + \frac{b\xi}{\gamma M_sL}(\mathbf{j}\cdot\nabla)\mathbf{m}.
\]
In particular, when an in-plane current is applied along the $x$-direction, the local field reduces to
\begin{equation}
\mathbf{f}(\mathbf{m}) = -Q(m_2\mathbf{e}_2 + m_3\mathbf{e}_3) + \mathbf{h}_e + \mathbf{h}_s
+ \frac{bJ}{\gamma M_sL}\mathbf{m}\times\mathbf{m}_x + \frac{bJ\xi}{\gamma M_sL}\mathbf{m}_x.
\label{equ:effectivedimensionless_withSTT}
\end{equation}

\section{The proposed method}
\label{sec:numericalmethods}

In this section, we shall describe the proposed method in details. The GSPM~\cite{panchi2020GSPM} is introduced firstly for completeness and comparison. For the LL equation \eqref{equ:LLdimensionless} with field \eqref{equ:effectivefield_withoutSTT} and Neumann boundary condition \eqref{equ:Neumannboundary}, the fractional step is applied
\begin{equation}
\begin{aligned}
\frac{\mathbf{m}^* - \mathbf{m}^n}{\Delta t} &= \epsilon\Delta_h\mathbf{m}^* + \mathbf{f}(\mathbf{m}^n),\\
\frac{\mathbf{m}^{n+1} - \mathbf{m}^n}{\Delta t} &= -\mathbf{m}^n\times\frac{\mathbf{m}^* - \mathbf{m}^n}{\Delta t} - \alpha\mathbf{m}^n\times\left(\mathbf{m}^n\times\frac{\mathbf{m}^* - \mathbf{m}^n}{\Delta t}\right),
\end{aligned}
\label{equ:fractionalstep}
\end{equation}
where $\Delta_h$ represents a discrete approximation to the Laplacian operator. In 3D, we use the second-order centered difference
\begin{equation}
\begin{aligned}
\Delta_h\mathbf{m}_{i,j,k} = &\frac{\mathbf{m}_{i+1,j,k}-
	2\mathbf{m}_{i,j,k}+\mathbf{m}_{i-1,j,k}}{\Delta x^2} +\\
&\frac{\mathbf{m}_{i,j+1,k}-2\mathbf{m}_{i,j,k}+
	\mathbf{m}_{i,j-1,k}}{\Delta y^2} +\\
&\frac{\mathbf{m}_{i,j,k+1}-2\mathbf{m}_{i,j,k}+
	\mathbf{m}_{i,j,k-1}}{\Delta z^2},
\end{aligned}
\label{LaplacianDiscrete}
\end{equation}
where $\mathbf{m}_{i,j,k}=\mathbf{m}((i - \frac 1{2})\Delta x,
(j - \frac 1{2})\Delta y, (k - \frac 1{2})\Delta z)$ with $i=0,1,\cdots,M,M+1$,
$j=0,1,\cdots,N,N+1$, and $k=0,1,\cdots,K,K+1$ being the indices of grid points in $x$-, $y$- and $z$-directions, respectively.
For the Neumann boundary condition \eqref{equ:Neumannboundary}, a second-order approximation yields
\begin{align*}
\mathbf{m}_{0,j,k} & =\mathbf{m}_{1,j,k},\quad \mathbf{m}_{M,j,k}  =\mathbf{m}_{M+1,j,k},\quad j = 1,\cdots,N,k=1,\cdots,K, \\
\mathbf{m}_{i,0,k} & =\mathbf{m}_{i,1,k},\quad \mathbf{m}_{i,N,k}  =\mathbf{m}_{i,N+1,k},\quad i = 1,\cdots,M,k=1,\cdots,K, \\
\mathbf{m}_{i,j,0} & =\mathbf{m}_{i,j,1},\quad \mathbf{m}_{i,j,K}  =\mathbf{m}_{i,j,K+1},\quad i = 1,\cdots,M,j=1,\cdots,N.
\end{align*}

Denote $\mathcal{L} = (I - \epsilon\Delta t\Delta_h)^{-1}$ with $I$ the identity operator. \eqref{equ:fractionalstep} can be rewritten as
\begin{gather}
m_i^* = \mathcal{L}(m_i^n + \Delta tf_i(\mathbf{m}^n)),\;\;i=1,2,3,\label{equ:fractionalstep1}\\
\label{equ:fractionalstep2}\mathbf{m}^{n+1} = \mathbf{m}^n - \mathbf{m}^n\times\mathbf{m}^* - \alpha\mathbf{m}^n\times\left(\mathbf{m}^n\times\mathbf{m}^*\right).
\end{gather}
Note that \eqref{equ:fractionalstep2} is updated in a Gauss-Seidel manner to avoid the stability issue \cite{NumGSPM2001}. It is nice that all linear systems in \eqref{equ:fractionalstep1} have constant coefficients with the symmetric and positive definite (spd) property, thus the optimal computational cost of solving linear systems here is $\mathcal{O}(n)$ where $n$ is the dofs. Due to the homogeneous Neumann boundary condition, the computational complexity is $\mathcal{O}(n\log(n))$ if discrete cosine transform is used to solve the linear systems here \cite{thesis2008Yanglei}. Meanwhile, FFT has also been used to calculate the stray field~\eqref{equ:scalarU} to reduce the computational complexity to $\mathcal{O}(n\log(n))$~\cite{CJreview2007}. Thus solving a linear system and updating the stray field have the same computational complexity in the asymptotic sense. However, the dofs of magnetization is actually $3n$. As a consequence,  updating the stray field is roughly three times costly compared to solving a linear system in \eqref{equ:fractionalstep1}.

\subsection{GSPM for LL equation}
In~\cite{panchi2020GSPM}, GSPM solves the LL equation in two steps:
\begin{itemize}
	\item Implicit Gauss-Seidel update
	\begin{equation}
	\begin{aligned}
	{g}_i^n &= \mathcal{L}({m}_i^n+\Delta tf_i(\mathbf{m}^n)),\quad i = 1,2,3, \\
	{g}_i^{*} &= \mathcal{L}({m}_i^{*}+\Delta tf_i(\mathbf{m}^*)),\quad i = 1,2,
	\end{aligned}
	\label{eqn:schemeA1}
	\end{equation}
	\begin{equation}
	\label{eqn:schemeA2}
	\begin{pmatrix}
	{m}_1^* \\{m}_2^* \\{m}_3^*
	\end{pmatrix} =
	\begin{pmatrix}
	m_1^n - (m_2^ng_3^n - m_3^ng_2^n) -
	\alpha(m_1^ng_1^n+m_2^ng_2^n + m_3^ng_3^n)m_1^n + \alpha g_1^n \\
	m_2^n - (m_3^ng_1^* - m_1^*g_3^n) -
	\alpha(m_1^{*}g_1^{*}+m_2^ng_2^n + m_3^ng_3^n)m_2^n + \alpha g_2^n \\
	m_3^n - (m_1^*g_2^* - m_2^*g_1^*) -
	\alpha(m_1^{*}g_1^{*}+m_2^{*}g_2^{*} + m_3^ng_3^n)m_3^n + \alpha g_3^n
	\end{pmatrix}.
	\end{equation}
	\item Projection onto $S^2$
	\begin{align}
	\label{eqn:schemeA3}
	\begin{pmatrix}
	{m}_1^{n+1} \\{m}_2^{n+1} \\{m}_3^{n+1}
	\end{pmatrix} = \frac 1{|\mathbf{m}^{*}|}
	\begin{pmatrix}
	{m}_1^{*} \\{m}_2^{*} \\{m}_3^{*}
	\end{pmatrix}.
	\end{align}
\end{itemize}
\begin{remark}
	In~\cite{seim2005}, the semi-implicit BDF1 scheme for the LL equation reads as
	\begin{equation}
	\frac{\mathbf{m}^{n+1} - \mathbf{m}^n}{\Delta t} = -\mathbf{m}^{n}\times(\epsilon\Delta\mathbf{m}^{n+1} +\mathbf{f}(\mathbf{m}^{n}))- \alpha\mathbf{m}^{n}\times(\mathbf{m}^{n}\times(\epsilon\Delta\mathbf{m}^{n+1} +\mathbf{f}(\mathbf{m}^{n}))).
	\label{equ:BDF1}
	\end{equation}
	A projection step is applied after \eqref{equ:BDF1} at each step. In fact, GSPM for the LL equation \eqref{eqn:schemeA1}-\eqref{eqn:schemeA3} can be obtained by applying the splitting strategy \eqref{equ:fractionalstep1}-\eqref{equ:fractionalstep2} to \eqref{equ:BDF1} and then updating \eqref{equ:fractionalstep2} in the Gauss-Seidel manner.
	This motivates the current work of applying the splitting strategy and updating in the Gauss-Seidel manner for the semi-implicit BDF2 scheme.
	\label{rk:BDF1statement}
\end{remark}
\begin{remark}
	From \eqref{eqn:schemeA1}, we know that at each time step GSPM solves linear systems of equations with $n$ dofs five times and updates the stray field with $3n$ dofs three times.
\end{remark}

\subsection{GSPM-BDF2 scheme for LL equation}

In~\cite{Changjian2019semiimplicit,jingrun2019analysis}, the semi-implicit BDF2 scheme for the LL equation reads as
\begin{multline}
\label{equ:semiBDF2}
\frac{\frac{3}{2}\mathbf{\hat{m}}^{n+2} - 2\mathbf{m}^{n+1} + \frac{1}{2}\mathbf{m}^n}{\Delta t} = -\mathbf{\tilde{m}}^{n+2}\times(\epsilon\Delta\mathbf{\hat{m}}^{n+2} +\mathbf{f}(\mathbf{\tilde{m}}^{n+2})), \\
- \alpha\mathbf{\tilde{m}}^{n+2}\times(\mathbf{\tilde{m}}^{n+2}\times(\epsilon\Delta\mathbf{\hat{m}}^{n+2} +\mathbf{f}(\mathbf{\tilde{m}}^{n+2})))
\end{multline}
\begin{equation}
\mathbf{m}^{n+2} = \frac{1}{\abs{\mathbf{\hat{m}}^{n+2}}}\mathbf{\hat{m}}^{n+2}
\end{equation}
with $\mathbf{\tilde{m}}^{n+2} = 2\mathbf{m}^{n+1} - \mathbf{m}^n$. This scheme has second-order accuracy in both space and time. At each step, it solves one linear system of equations with variable coefficients, a nonsymmetric structure, and $3n$ dofs. Note that only one update is needed for the stray field per step in \eqref{equ:semiBDF2}.

Similarly, we solve the semi-implicit BDF2 scheme \eqref{equ:semiBDF2} using the splitting strategy and updating in the Gauss-Seidel manner.
The Gauss-Seidel update is applied for the following equation after splitting
\begin{equation}
\frac{\frac{3}{2}\mathbf{\hat{m}}^{n+2} - 2\mathbf{m}^{n+1} + \frac{1}{2}\mathbf{m}^n}{\Delta t} = -\mathbf{\tilde{m}}^{n+2}\times\mathcal{L}\mathbf{\tilde{m}}^{n+2} - \alpha\mathbf{\tilde{m}}^{n+2}\times(\mathbf{\tilde{m}}^{n+2}\times\mathcal{L}\mathbf{\tilde{m}}^{n+2}).
\label{equ:GSPM-BDF2}
\end{equation}
Therefore, in details, GSPM-BDF2 works as follows
\begin{gather*}
g_i^{n+2} = \mathcal{L}(\tilde{m}_i^{n+2} + \Delta tf_i(\mathbf{\tilde{m}}^{n+2})),\;\; i = 1, 2, 3,\\
\frac{3}{2}m_1^{*} = 2m_1^{n+1} - \frac{1}{2}m_1^n - (\tilde{m}^{n+2}_2g^{n+2}_3 - \tilde{m}_3^{n+2}g_2^{n+2}) - \alpha(\tilde{m}_1^{n+2}g_1^{n+2} + \tilde{m}_2^{n+2}g_2^{n+2} + \tilde{m}_3^{n+2}g_3^{n+2})\tilde{m}_1^{n+2} + \alpha g_1^{n+2},\\
\tilde{m}^{*}_1 = 2m_1^{*} - m_1^{n+1},\ \ g_1^{*} = \mathcal{L}(\tilde{m}^{*}_1 + \Delta tf_1(\mathbf{\tilde{m}}^{n+2})),\\
\frac{3}{2}m_2^{*} = 2m_2^{n+1} - \frac{1}{2}m_2^n - (\tilde{m}^{n+2}_3g^{*}_1 - \tilde{m}_1^{*}g_3^{n+2}) - \alpha(\tilde{m}_1^{*}g_1^{*} + \tilde{m}_2^{n+2}g_2^{n+2} + \tilde{m}_3^{n+2}g_3^{n+2})\tilde{m}_2^{n+2} + \alpha g_2^{n+2},\\
\tilde{m}^{*}_2 = 2m_2^{*} - m_2^{n+1},\ \ g_2^{*} = \mathcal{L}(\tilde{m}^{*}_2 + \Delta tf_2(\mathbf{\tilde{m}}^{n+2})),\\
\frac{3}{2}m_3^{*} = 2m_3^{n+1} - \frac{1}{2}m_3^n - (\tilde{m}^{*}_1g^{*}_2 - \tilde{m}_2^{*}g_1^{*}) - \alpha(\tilde{m}_1^{*}g_1^{*} + \tilde{m}_2^{*}g_2^{*} + \tilde{m}_3^{n+2}g_3^{n+2})\tilde{m}_3^{n+2} + \alpha g_3^{n+2},\\
\begin{pmatrix}
m_1^{n+2} \\ m_2^{n+2} \\ m_3^{n+2}
\end{pmatrix} = \frac{1}{\abs{\mathbf{m}^*}}
\begin{pmatrix}
m_1^{*} \\ m_2^{*} \\ m_3^{*}
\end{pmatrix}.
\end{gather*}

In \eqref{equ:semiBDF2}, one linear system with variable coefficients, a nonsymmetric structure, and $3n$ dofs is solved at each step. This is avoided in GSPM-BDF2 thanks to the splitting strategy. Five linear systems with constant coefficients, spd structure, and $n$ dofs are solved. Meanwhile, the Gauss-Seidel update is tested to be unconditionally stable. GSPM-BDF2 inherits the advantage of BDF2 that only one update of the stray field is needed with $3n$ dofs. We summarize the main computational costs of GSPM and GSPM-BDF2 in \cref{tab:costcomparsion} by counting the number of linear systems of equations and the number of updates for the stray field per step. Simulations suggest that one update of stray field with $3n$ dofs is computationally comparable to solving three linear systems of equations with $n$ dofs, thus GSPM-BDF2 saves about $60\%$ computational time over GSPM; see \Cref{sec:numericalsimulation} for details.
\begin{table}[htbp]
	\centering
	\caption{Main computational costs of GSPM and GSPM-BDF2 in micromagnetic simulations.}
	\begin{tabular}{||c|c|c||}
		\hline
		 & \#(linear systems of equations) (dofs) & \#(stray field updates) (dofs)\\
		 \hline
		 GSPM & 5 ($n$) & 3 ($3n$) \\
		 GSPM-BDF2 & 5 ($n$) & 1 ($3n$)\\
		 \hline
	\end{tabular}
\label{tab:costcomparsion}
\end{table}

GSPM-BDF2 is expected to be first-order accurate in time and second-order accurate in space due to the splitting strategy and the Gauss-Seidel update, both shall be verified numerically later. Compared to the semi-implicit BDF2 scheme, GSPM-BDF2 only solves linear systems of equations with constant coefficients and a spd structure. Many efficient linear solvers can be implemented directly, while the linear system in BDF2 has to be solved by GMRES which is an iterative solver and is not so efficient as solvers for spd matrices.

\begin{remark}
There are two improved GSPMs proposed in~\cite{panchi2020GSPM}: Scheme A and Scheme B. The current work proposes a more efficient method by combining Scheme A and BDF2. Scheme B solves three linear systems of equations with $n$ dofs and updates the stray field with $3n$ dofs three times. One may wonder how the combination of Scheme B and BDF2 works. Unfortunately, we find that the stability of the resulting scheme depends on the damping parameter $\alpha$ for micromagnetic simulations. The underlying reason is that Scheme B updates the magnetization in one step and does a projection step in a subsequent step and a delay of the update for the stray field leads to the instability with respect to the damping parameter.
\end{remark}

\section{Numerical simulations}
\label{sec:numericalsimulation}

\subsection{Accuracy check}
In 1D, we consider the following LL equation
\begin{equation}
\mathbf{m}_t = -\mathbf{m}\times\mathbf{m}_{xx} - \alpha\mathbf{m}\times(\mathbf{m}\times\mathbf{m}_{xx}).
\label{equ:LLforaccuracy}
\end{equation}
Choose the exact solution
\[
\mathbf{m}_e = (\cos(\bar{x})\sin(t), \sin(\bar{x})\sin(t), \cos(t))
\]
with $\bar{x} = x^2(1-x)^2$. A forcing term will be added on the right-hand side of \eqref{equ:LLforaccuracy} with $\mathbf{\hat{f}} = \mathbf{m}_{et} + \mathbf{m}_e\times\mathbf{m}_{exx} + \alpha\mathbf{m}_e\times(\mathbf{m}_e\times\mathbf{m}_{exx})$. The convergence rate with respect to the temporal step size and the spatial step size is recorded in \cref{tab:1daccuracy}.
\begin{table}[htbp]
	\centering
	\caption{Convergence rates in time and space for the 1D example. The final time $T=1.0e-02$ and the damping parameter $\alpha = 0.01$. $\Delta x = 1.0e-03$ is used for the temporal accuracy and $\Delta t = 1.0e-06$ is used for the spatial accuracy, respectively.}
	\label{tab:1daccuracy}
	\begin{tabular}{||c|c|cccc|c||}
		\hline
		\multirow{3}{*}{\makecell[c]{Temporal\\accuracy}} & nt & 1000 & 2000 & 4000 & 8000 & order\\
		\cline{2-7}
		& GSPM & 6.01e-08 & 3.01e-08 & 1.52e-08 & 7.72e-09 & 0.99\\
		& GSPM-BDF2 & 6.01e-08 & 3.02e-08 & 1.52e-08 & 7.73e-09 & 0.99\\
		\cline{1-7}
		\multirow{3}{*}{\makecell[c]{Spatial\\accuracy}} & nx & 10 & 20 & 40 & 80 & order\\
		\cline{2-7}
		& GSPM & 1.24e-06 & 4.24e-07 & 1.27e-07 & 4.03e-08 & 1.66\\
		& GSPM-BDF2 & 1.24e-06 & 4.24e-07 & 1.27e-07 & 4.03e-08 & 1.66\\
		\hline
	\end{tabular}
\end{table}

In 3D, we choose the exact solution
\[
\mathbf{m}_e = (\cos(\bar{x}\bar{y}\bar{z})\sin(t), \sin(\bar{x}\bar{y}\bar{z})\sin(t), \cos(t)),
\]
where $\bar{x} = x^2(1-x)^2$, $\bar{y} = y^2(1-y)^2$ and $\bar{z} = z^2(1-z)^2$. Uniform discretization over $\Omega = [0, 1]^3$ is used for the spatial accuracy. $\Omega=[0, 2]\times[0, 1]\times[0, 0.2]$ with $128\times64\times10$ cubes is used for the temporal accuracy. Errors are recorded in \cref{tab:3daccuracy}.
From 1D and 3D results, the accuracy of GSPM-BDF2 is $\mathcal{O}(\Delta t + (\Delta x)^2)$. Moreover, the scheme is tested to be unconditionally stable by varying the temporal step size.
\begin{table}[htbp]
	\centering
	\caption{Convergence rates in time and space for the 3D example. The final time $T=1.0e-05$ and the damping parameter $\alpha = 0.01$. $\Delta x = 1.0e-03$ is used for the temporal accuracy and $\Delta t = 1.0e-09$ is used for the spatial accuracy, respectively.}
	\label{tab:3daccuracy}
	\begin{tabular}{||c|c|cccc|c||}
		\hline
		\multirow{3}{*}{\makecell[c]{Temporal\\accuracy}} & nt & 10 & 20 & 40 & 80 & order\\
		\cline{2-7}
		& GSPM & 1.00e-06 & 5.00e-07 & 2.50e-07 & 1.25e-07 & 1.00\\
		& GSPM-BDF2 & 1.00e-06 & 5.00e-07 & 2.50e-07 & 1.25e-07 & 1.00\\
		\cline{1-7}
		\multirow{3}{*}{\makecell[c]{Spatial\\accuracy}} & nx=ny=nz & 6 & 8 & 10 & 12 & order\\
		\cline{2-7}
		& GSPM & 2.91e-14 & 1.72e-14 & 1.13e-14 & 7.92e-15 & 1.88\\
		& GSPM-BDF2 & 2.91e-14 & 1.72e-14 & 1.13e-14 & 7.92e-15 & 1.88\\
		\hline
	\end{tabular}
\end{table}

\subsection{Stability with respect to the damping parameter}

The original GSPM is found to be conditionally stable with respect to the damping parameter~\cite{ImproveGSPM2003} if the stray field is updated only once per step~\cite{NumGSPM2001}. Here we examine the performance of GSPM-BDF2 for small damping parameters. Consider the magnetization dynamics of a ferromagnet over $\Omega=1\;\mu\textrm{m}\times1\;\mu\textrm{m}\times0.02\;\mu\textrm{m}$ and the final time $1.6$ nanoseconds. The spatial mesh size is $4\;\textrm{nm}\times4\;\textrm{nm}\times4\;\textrm{nm}$, and the temporal step size is $1$ picosecond.

\begin{figure}[htbp]
	\centering
	\subfloat[Angle profile]{\includegraphics[width=2.5in]{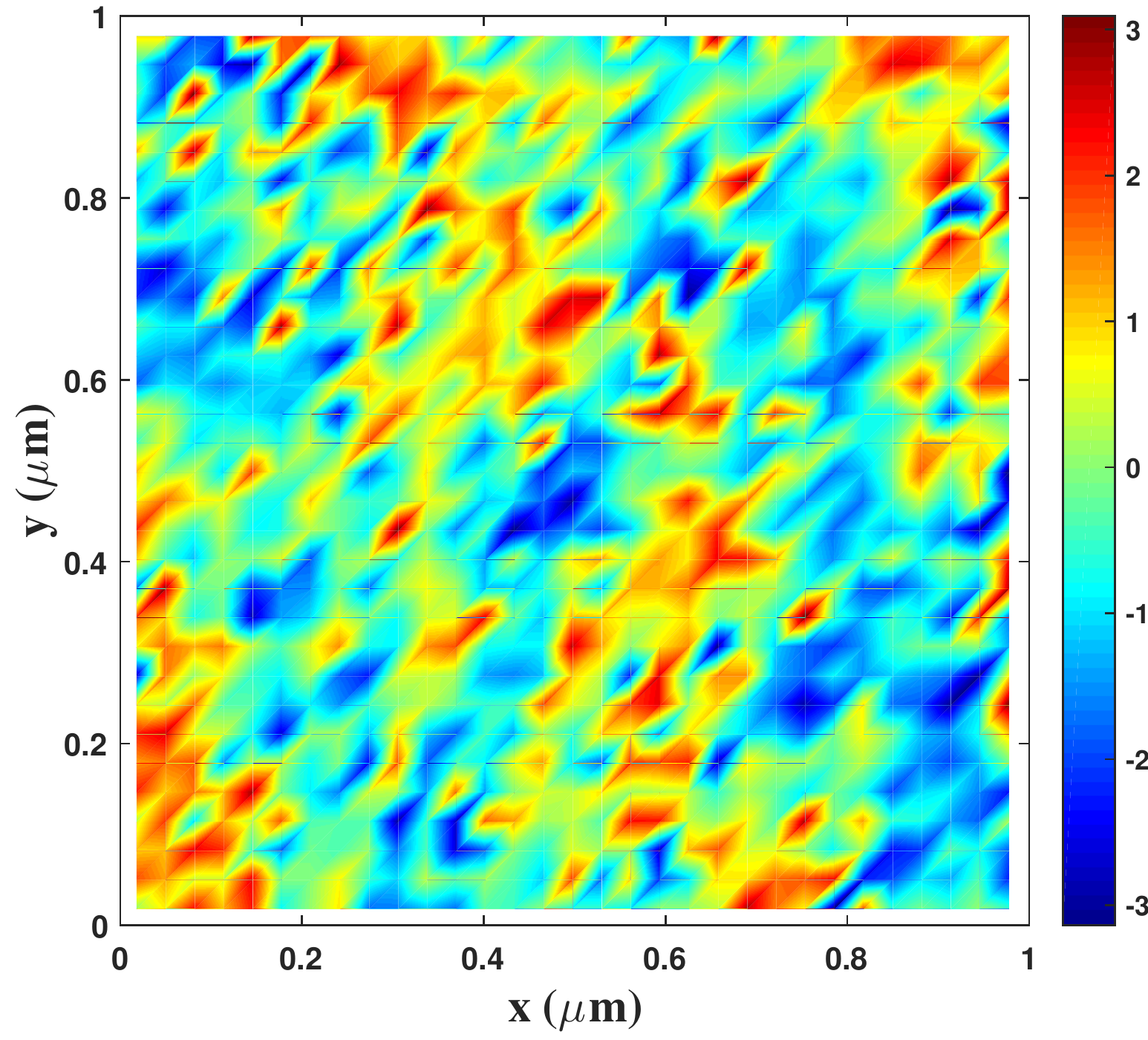}}
	\subfloat[Magnetization profile]{\includegraphics[width=2.25in]{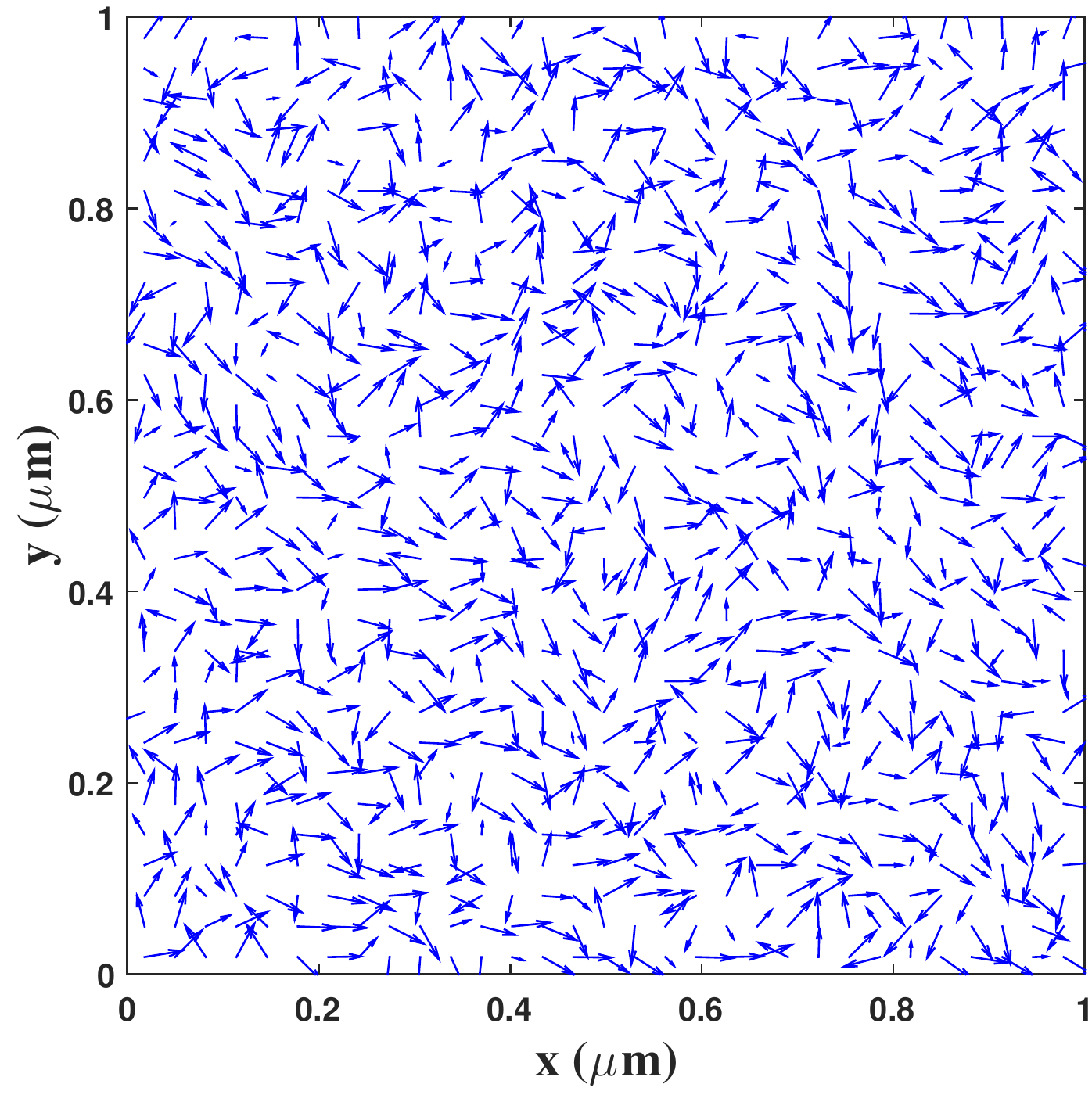}}
	\caption{Simulation results of the LL equation using GSPM with only one update of stray field at each step. The magnetization on the centered slice of the material in the $xy$ plane is used. Left: a color plot of the angle between the in-plane magnetization and the $x$ axis; Right : an arrow plot of the in-plane magnetization.}
	\label{fig:incorrectdynamics}
\end{figure}

\begin{figure}[htbp]
	\centering
	\subfloat[Angle profile]{\includegraphics[width=2.5in]{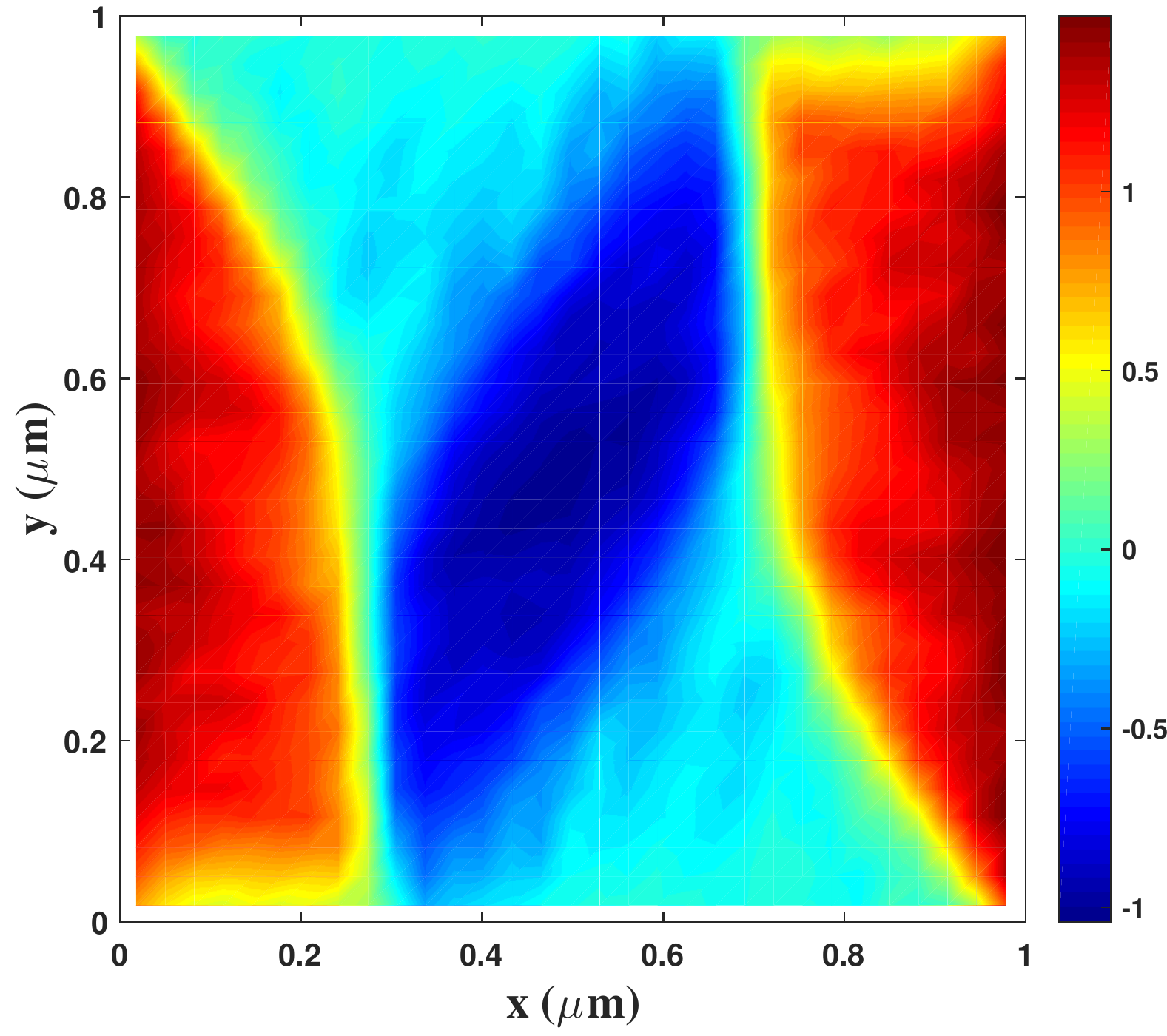}}
	\subfloat[Magnetization profile]{\includegraphics[width=2.25in]{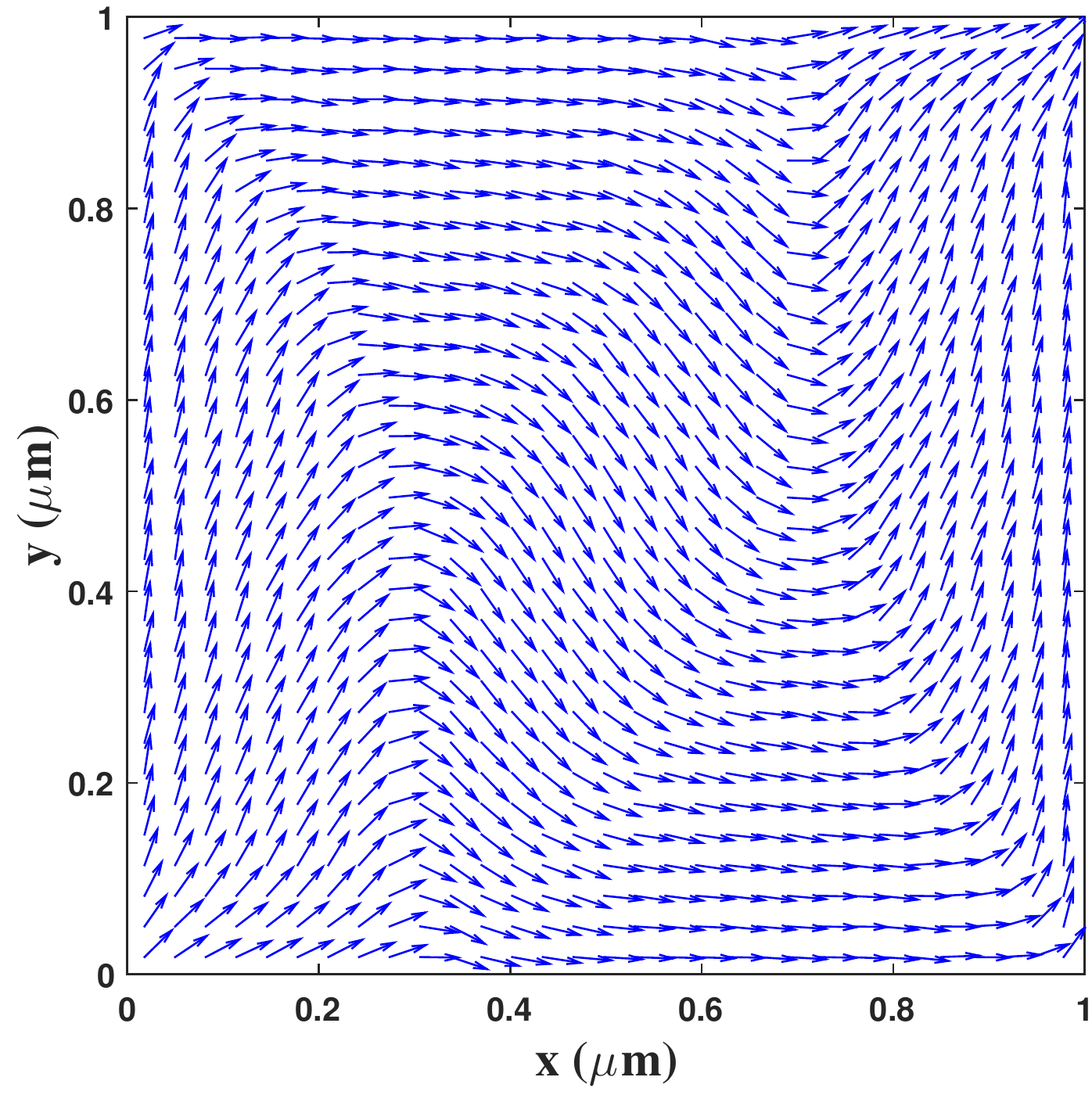}}
	\caption{Simulation results of the LL equation using GSPM with three updates of stray field at each step. The magnetization on the centered slice of the material in the $xy$ plane is used. Left: a color plot of the angle between the in-plane magnetization and the $x$ axis; Right : an arrow plot of the in-plane magnetization.}
	\label{fig:correctdynamicsSchemeA}
\end{figure}

For comparison, we first update the stray field in GSPM \eqref{eqn:schemeA1}-\eqref{eqn:schemeA3} only once, i.e., only $\mathbf{f}(\mathbf{m}^n)$ is used per step.
Results are shown in \Cref{fig:incorrectdynamics} when $\alpha = 0.01$. The whole simulation takes $493.75$ seconds. If we use the GSPM \eqref{eqn:schemeA1}-\eqref{eqn:schemeA3} with three updates of the stray field, we obtain the simulation results in \Cref{fig:correctdynamicsSchemeA}. The whole simulation takes $1225.92$ seconds. Comparing these two results, we find that updating the stray field only once leads to the numerical instability with respect to the damping parameter in micromagnetic simulations, although the computational cost is reduced by $60\%$.

\begin{figure}[htbp]
	\centering
	\subfloat[Angle profile]{\includegraphics[width=2.5in]{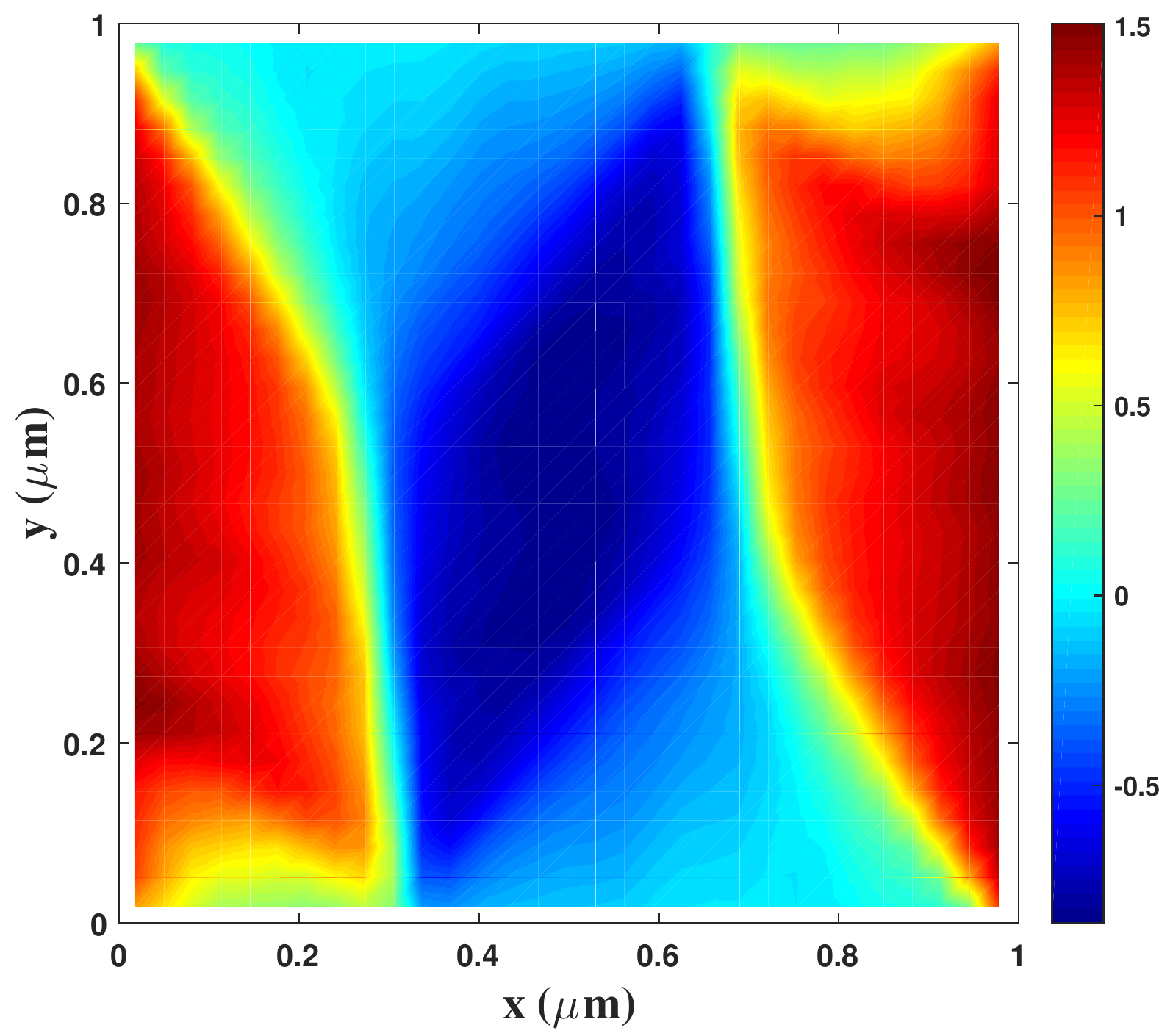}}
	\subfloat[Magnetization profile]{\includegraphics[width=2.25in]{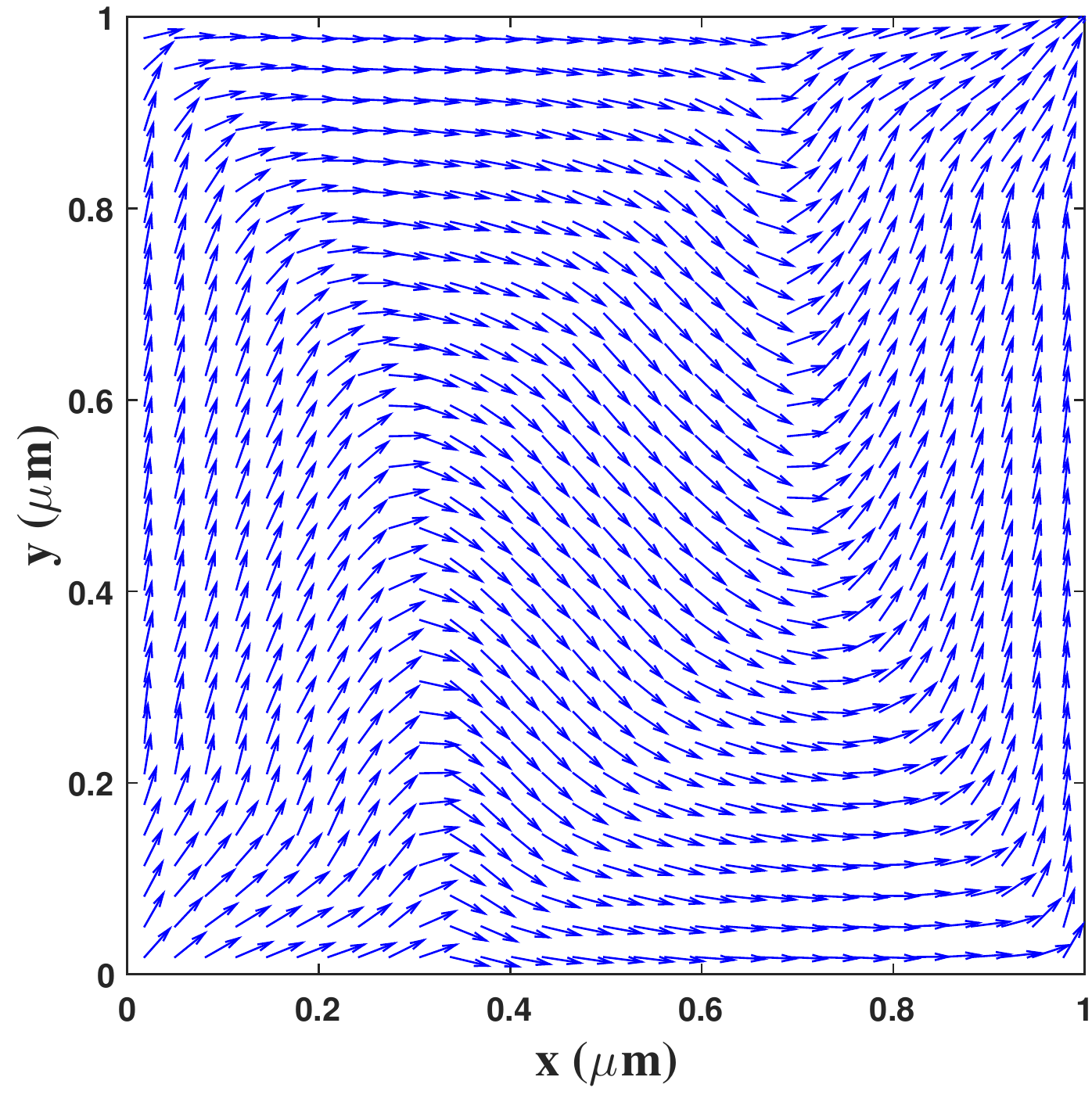}}
	\caption{Simulation results of the LL equation using GSPM-BDF2. The magnetization on the centered slice of the material in the $xy$ plane is used. Left: a color plot of the angle between the in-plane magnetization and the $x$ axis; Right : an arrow plot of the in-plane magnetization.}
	\label{fig:GSPM-BDF2sstate}
\end{figure}

Simulation results of GSPM-BDF2 are plotted in \cref{fig:GSPM-BDF2sstate} when $\alpha = 0.01$. It is clear that GSPM-BDF2 is able to produce correct magnetization dynamics though the stray field is updated only once per step. The whole simulation takes $511.39$ seconds and saves $58\%$ computational time over GSPM \eqref{eqn:schemeA1}-\eqref{eqn:schemeA3}. This study shows that the proposed method is stable with respect to the damping parameter even the number of updates for the stray field is minimized.

\subsection{Standard Problem \#4}
\begin{figure}[htbp]
	\centering
	\includegraphics[width=5.in]{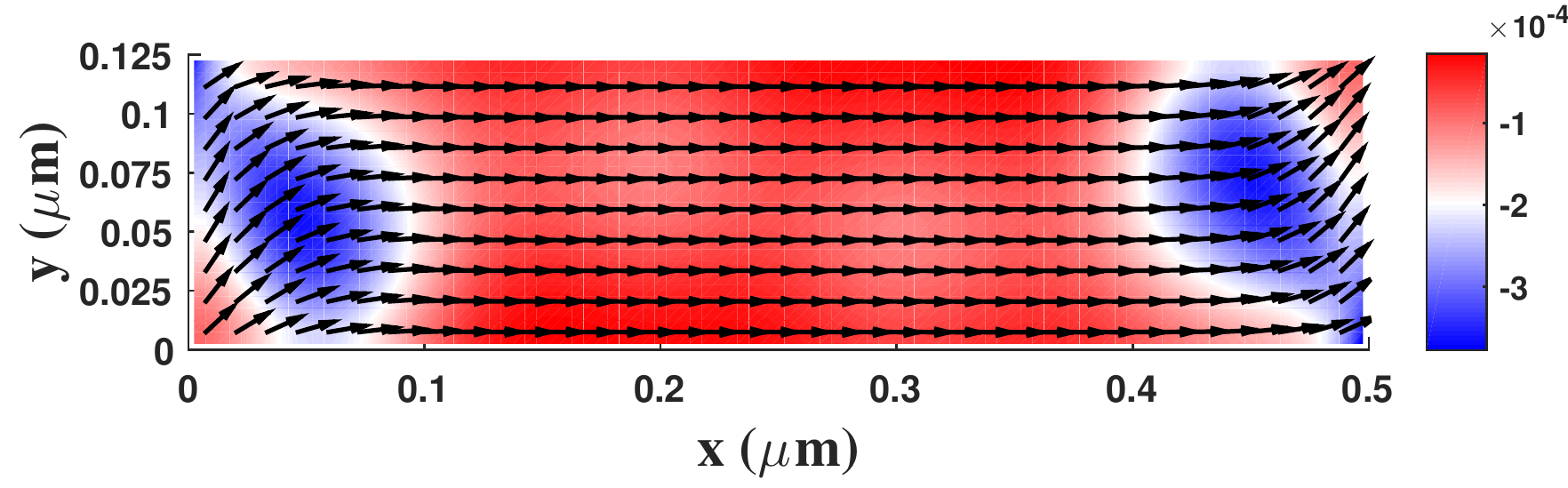}
	\caption{Magnetization profile of the initial s-state. This is generated for a random initial state under a magnetic field $100\;\mathrm{mT}$ along the $[1,1,1]$ direction with a successive reduction to $0$.}\label{fig:initial s state}
\end{figure}

The setup of this benchmark problem is as follows: film size $500\;\mathrm{nm}\times125\;\mathrm{nm}\times3\;\mathrm{nm}$;  initial state: an equilibrium s-state.
A magnetic field with sufficient magnitude is applied to reverse the magnetization of the s-state. Following the description of Standard Problem \#4, we first generate the initial s-state in \cref{fig:initial s state}.

Two different magnetic fields are applied:
\begin{itemize}
	\item A field 25 $\mathrm{mT}$, directed 170 degrees counterclockwise from the positive $\mathrm{x}$ axis;
	\item A field 36 $\mathrm{mT}$, directed 190 degrees counterclockwise from the positive $\mathrm{x}$ axis.
\end{itemize}
Two different mesh strategies are used: a coarse mesh with cell size $5\;\mathrm{nm}\times5\;\mathrm{nm}\times3\;\mathrm{nm}$ and a fine mesh with cell size $2.5\;\mathrm{nm}\times2.5\;\mathrm{nm}\times3\;\mathrm{nm}$. All the results shown here are independent of the mesh strategy.

\begin{figure}[htbp]
	\centering
	\subfloat[Magnetization dynamics]{\includegraphics[width=2.6in]{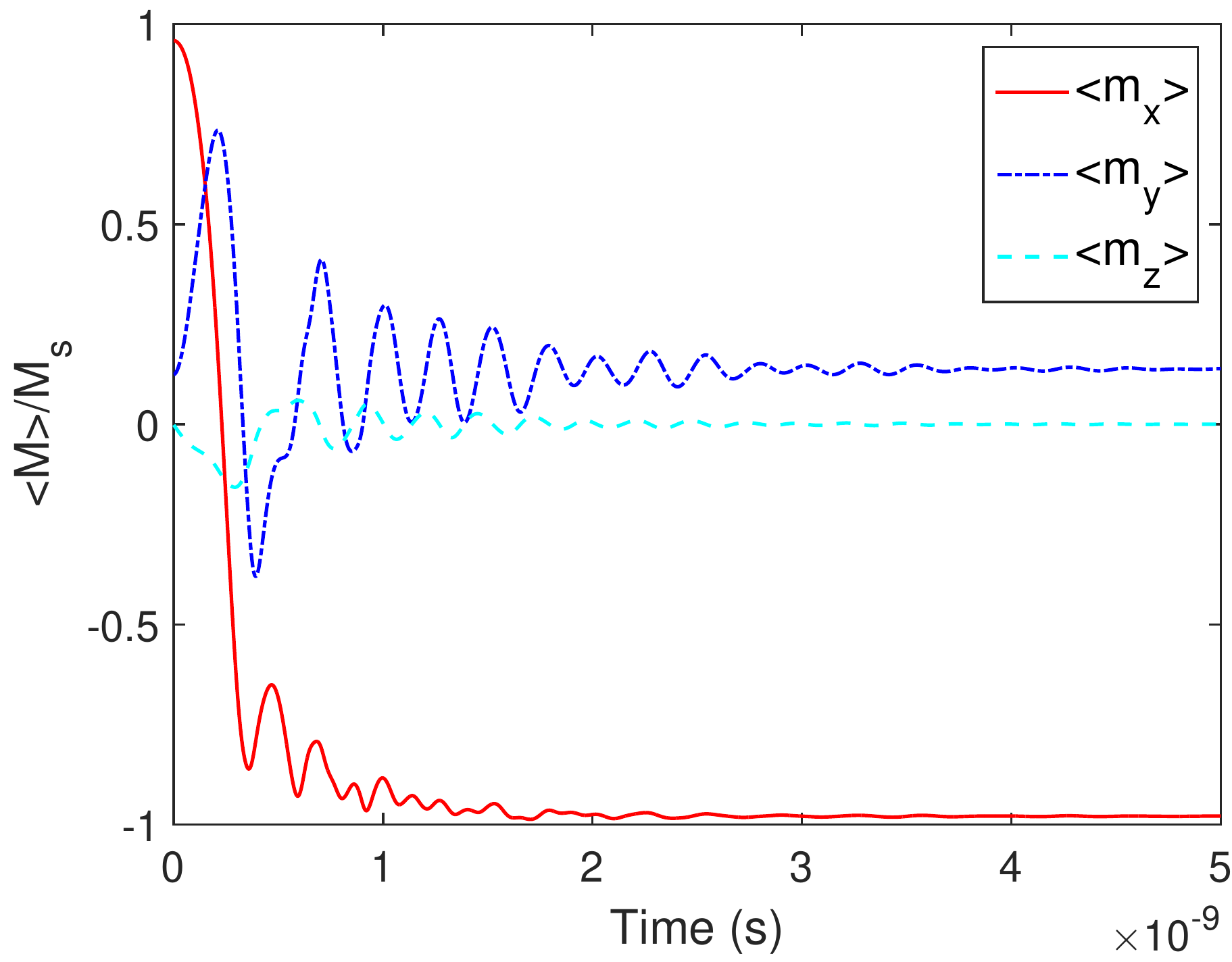}}
	\subfloat[Comparison of $<m_y>$]{\label{subfig:myfield1}\includegraphics[width=2.6in]{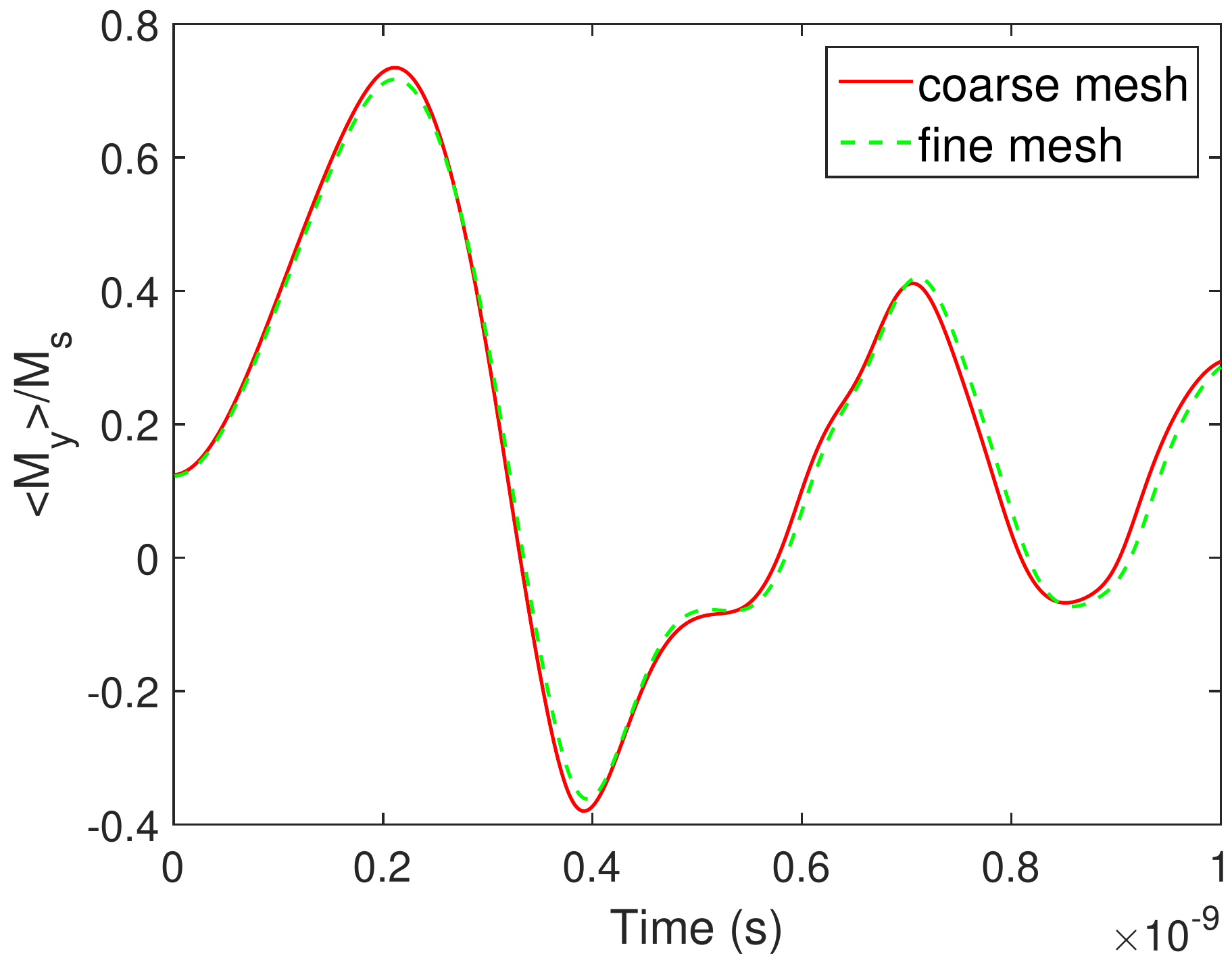}}
	\caption{Left: dynamics of the spatially averaged magnetization on the coarse mesh under the external field of $25\;\mathrm{mT}$; Right: comparison of the averaged $<m_y>$ on two different meshes.}
	\label{fig:std4field1}
\end{figure}

\begin{figure}[htbp]
	\centering
	\includegraphics[width=5.5in]{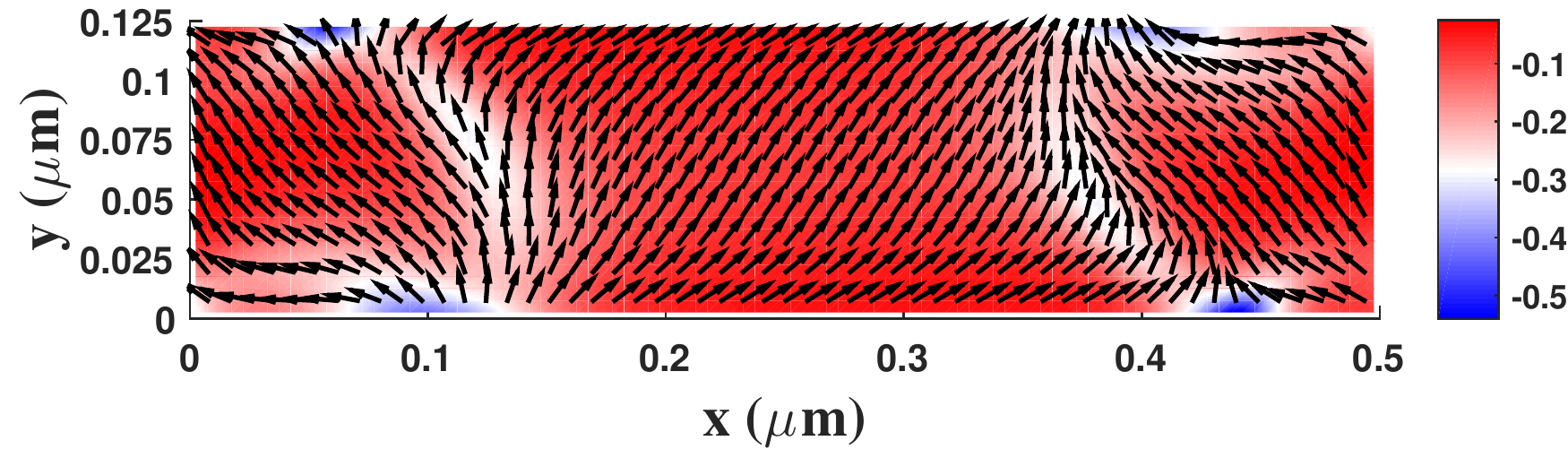}
	\caption{Magnetization profile when the averaged $<m_x>=0$ under the external field of $25\;\mathrm{mT}$. The color map is given by the $\mathrm{z}$-component of the magnetization.}
	\label{fig:field1_reveal_phase}
\end{figure}

Dynamics of the spatially averaged magnetization under the external field of $25\;\mathrm{mT}$ is plotted in \cref{fig:std4field1}.
Magnetization profile when the averaged $<m_x>=0$ under the external field of $25\;\mathrm{mT}$ is visualized in \cref{fig:field1_reveal_phase}. The color map is given by the $\mathrm{z}$-component of the magnetization.

\begin{figure}[htbp]
	\centering
	\subfloat[Magnetization dynamics]{\includegraphics[width=2.6in]{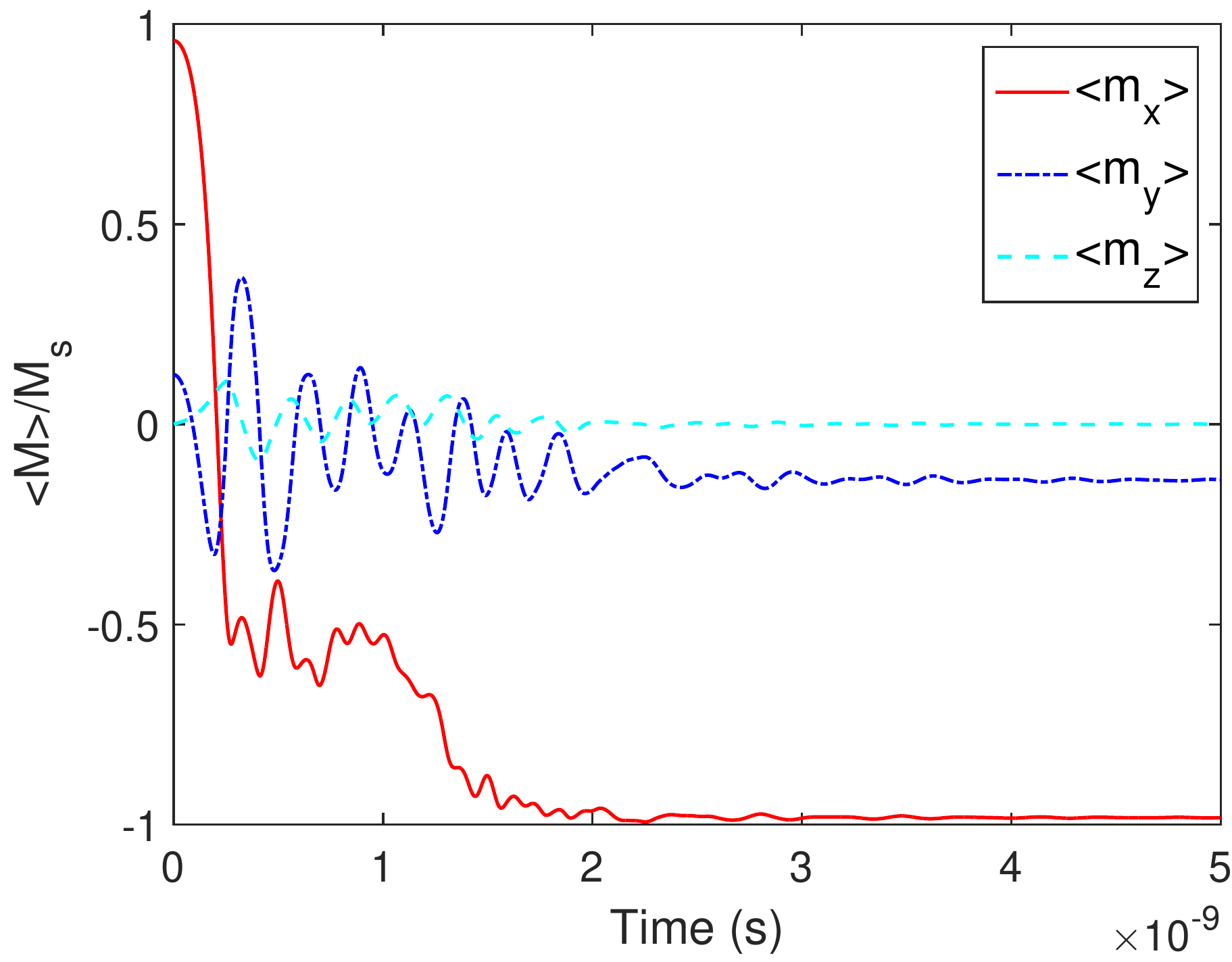}}
	\subfloat[Comparison of $<m_y>$]{\label{subfig:field2}\includegraphics[width=2.6in]{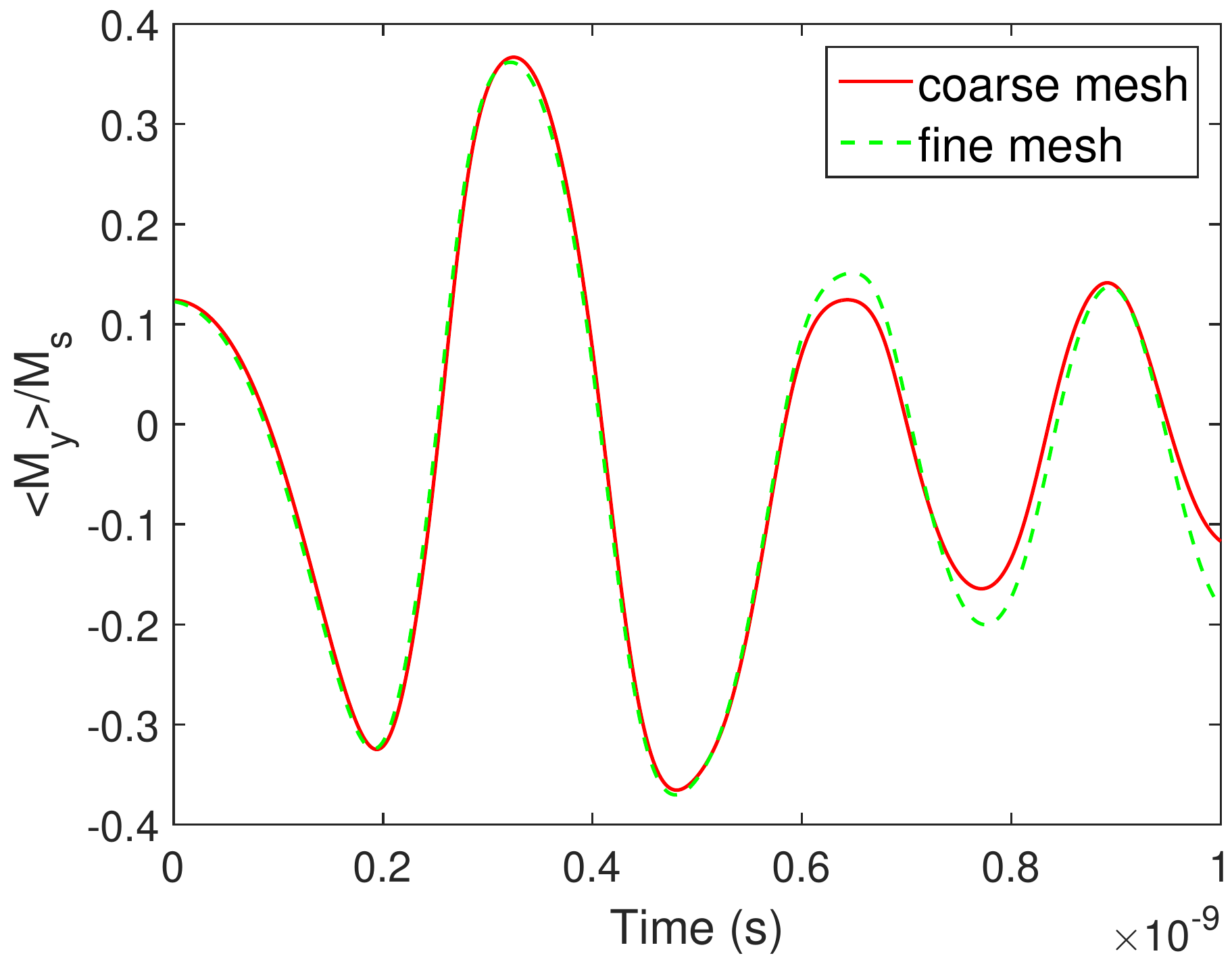}}
	\caption{Left: dynamics of the spatially averaged magnetization on the coarse mesh under the external field of $36\;\mathrm{mT}$; Right: comparison of the averaged $<m_y>$ on two different meshes.}
	\label{fig:std4field2}
\end{figure}

\begin{figure}
	\centering
	\includegraphics[width=5.5in]{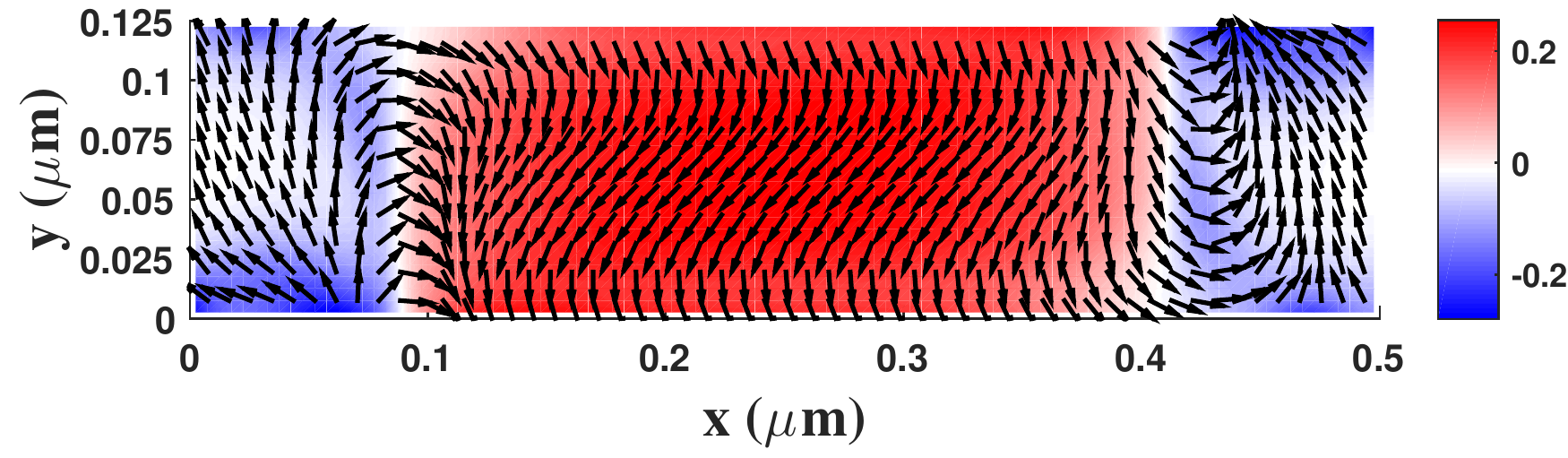}
	\caption{Magnetization profile when the averaged $<m_x>=0$ under the external field of $36\;\mathrm{mT}$. The color map is given by the $\mathrm{z}$-component of the magnetization.}
	\label{fig:field2_reveal_phase}
\end{figure}

In the presence of the $36\;\mathrm{mT}$ field, dynamics of the spatially averaged magnetization is plotted in \cref{fig:std4field2} and the magnetization profile when the averaged $<m_x>=0$ is visualized in \cref{fig:field2_reveal_phase}, respectively.

When the field is applied at 170 degrees, the magnetization at the center of the rectangle rotates in the same direction as that at the two ends during the reversal process. When the field is applied at 190 degrees, the magnetization at the center of the rectangle rotates in the opposite direction as that at the two ends during the reversal process. These are in good agreements with the reports listed in~\cite{NIST}. To check the efficiency, we record the computational costs of the proposed method and OOMMF in \cref{tab:ex4}. The current work saves $82\%$ computational time over OOMMF.
\begin{table}[htbp]
	\centering
	\caption{Computational costs of the proposed method and OOMMF for Standard Problem \#4~(unit: seconds) when the coarse mesh is used.}\label{tab:ex4}
	\begin{tabular}{||c|cc|c||}
		\hline
	Standard Problem \#4	& The proposed method & OOMMF  & Saving\\
		\hline
		field $25\;\mathrm{mT}$ & 20.47 & 115.32 & 82\% \\
		field $36\;\mathrm{mT}$ & 20.33 & 116.41 & 83\% \\
		\hline
	\end{tabular}
\end{table}

\subsection{Standard Problem \#5}
Standard Problem \#5 considers the magnetization dynamics in the presence of STT~\cite{ZhangLi2004STT}, which is modeled by the LL equation with \eqref{equ:effectivedimensionless_withSTT}. Since Neumann boundary condition is used, zero-spin torques of the magnetization on the boundaries where the spin current enters and leaves is naturally satisfied. The setup is as follows: ferromagnet size $100\;\mathrm{nm}\times100\;\mathrm{nm}\times10\;\mathrm{nm}$; cell size $2\mathrm{nm}\;\times2\mathrm{nm}\;\times2\;\mathrm{nm}$. At each point $(x, y, z)\in\Omega$, the initial magnetization is chosen as $\mathbf{m} = \mathbf{g}/\abs{\mathbf{g}}$, where
\begin{equation}
\mathbf{g}(x,y,z) = [-y, x, R]
\label{equ:std5ini}
\end{equation}
and $R = 10\;\mathrm{nm}$. The equilibrium configuration after relaxation is chosen as the initial state (vortex pattern) for all simulations in what follows. 

\begin{figure}[htbp]
	\centering
	\subfloat[Magnetization dynamics]{\includegraphics[width=2.8in]{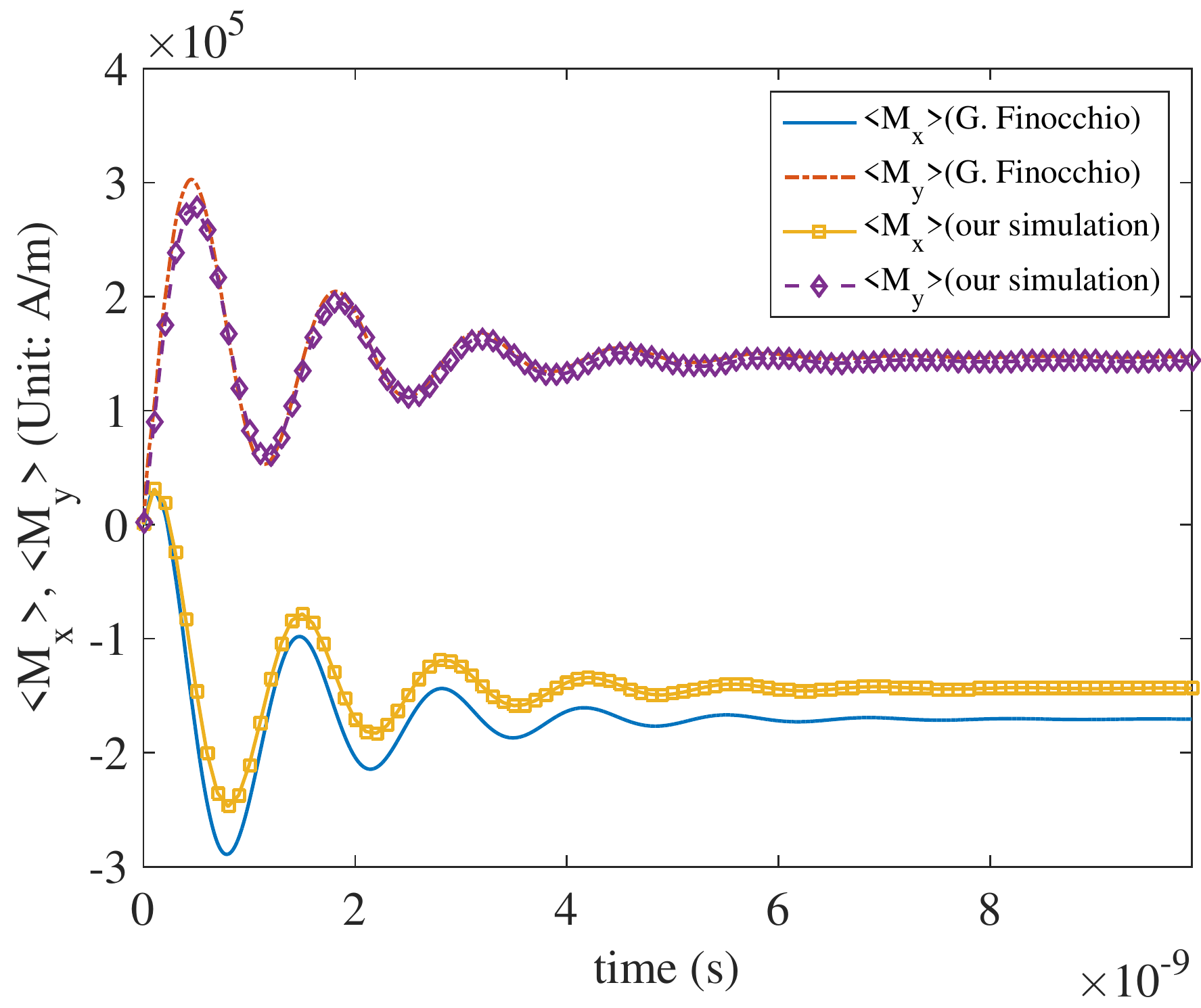}}
	\subfloat[Final state]{\includegraphics[width=2.6in]{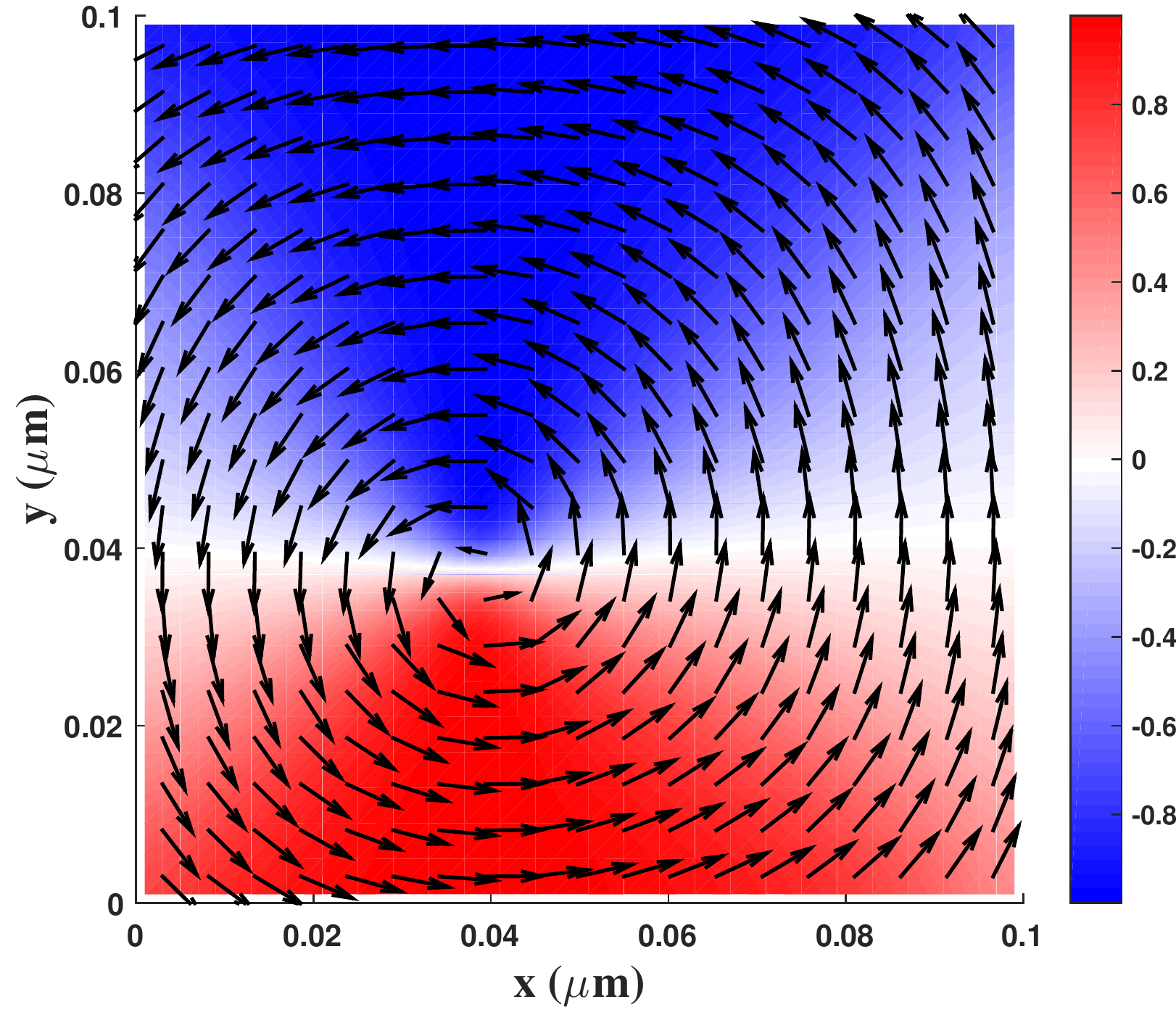}}
	\caption{Magnetization dynamics and the final state in the case of $\xi=0.5$ and $bJ=72.45 \;\mathrm{m}/\mathrm{s}$. Left: dynamics of the spatially averaged magnetization with the result of G. Finocchio et al. as the reference; Right: the final state colored by the $x$-component of magnetization. The first component of the spatially averaged magnetization $<M_x>$ at $10\;\mathrm{ns}$ is $-1.43\times10^5\;\mathrm{A}/\mathrm{m}$.}
	\label{fig:xi=05andu=72.45}
\end{figure}

There are four sets of parameters used in Standard Problem \#5:
\begin{itemize}
	\item[1)] $bJ = 72.35 \mathrm{m}/\mathrm{s}$ and $\xi = 0$;
	\begin{figure}[htbp]
		\centering
		\subfloat[Magnetization dynamics]{\includegraphics[width=2.8in]{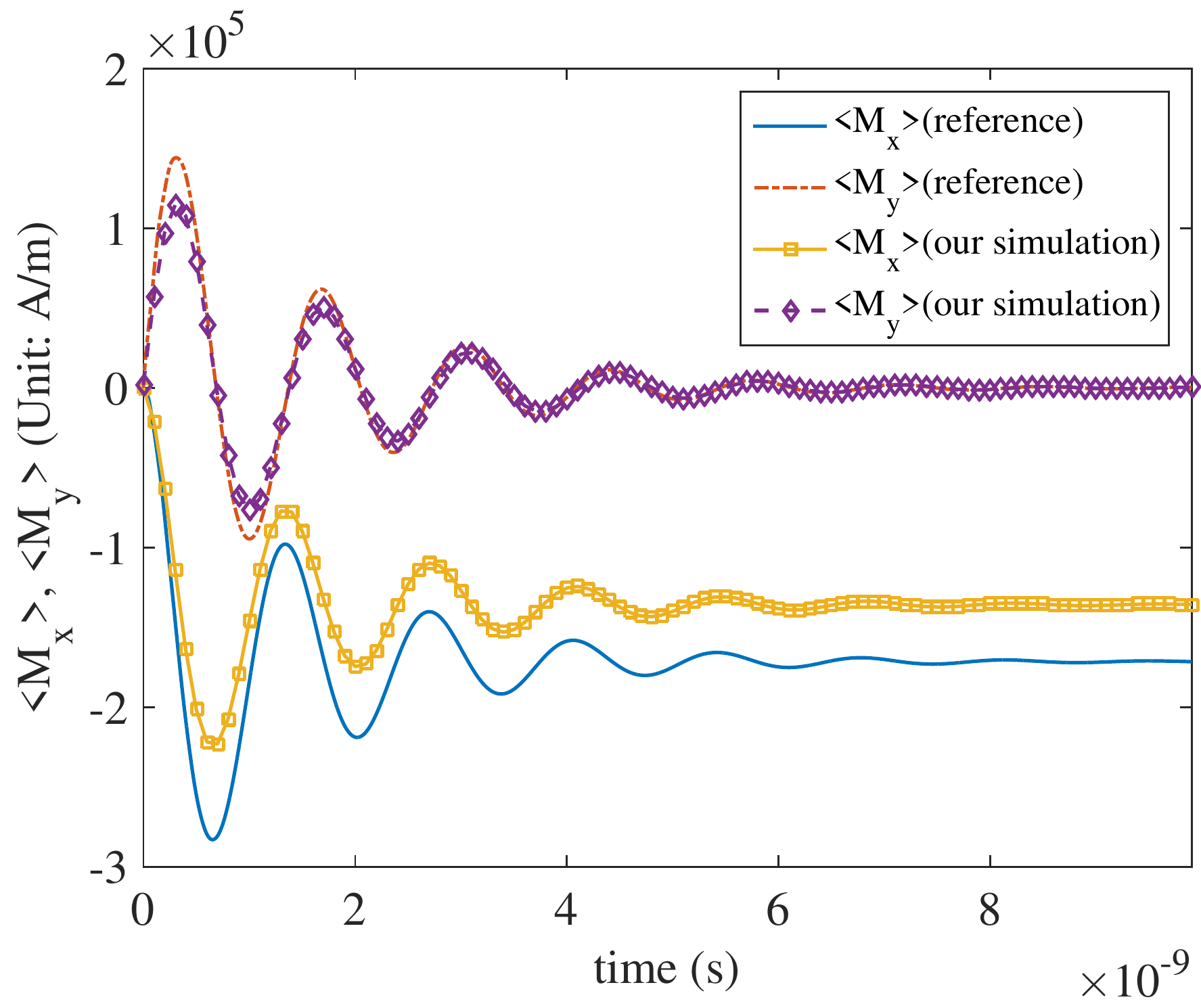}}
		\subfloat[Final state]{\includegraphics[width=2.6in]{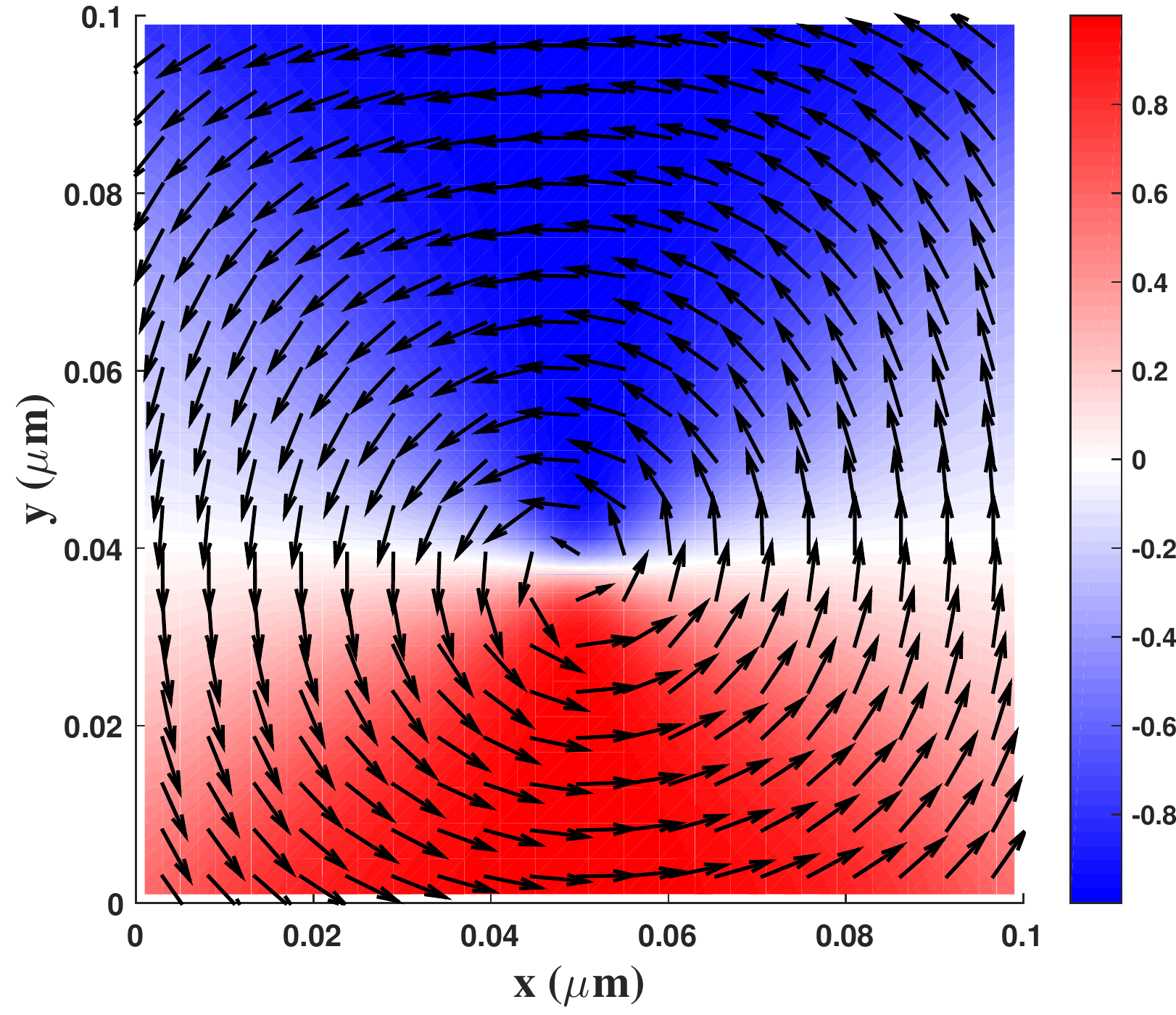}}
		\caption{Magnetization dynamics and the final state in the case of $\xi = 0$ and $bJ = 72.35 \;\mathrm{m}/\mathrm{s}$.
			Left: dynamics of the spatially averaged magnetization with the result of D. G. Porter as the reference; Right: the final state colored by the $x$-component of magnetization.}
		\label{fig:xi=0}
	\end{figure}
	
	\item[2)] $bJ = 72.17 \mathrm{m}/\mathrm{s}$ and $\xi = 0.05$;
	\begin{figure}[htbp]
		\centering
		\subfloat[Magnetization dynamics]{\includegraphics[width=2.8in]{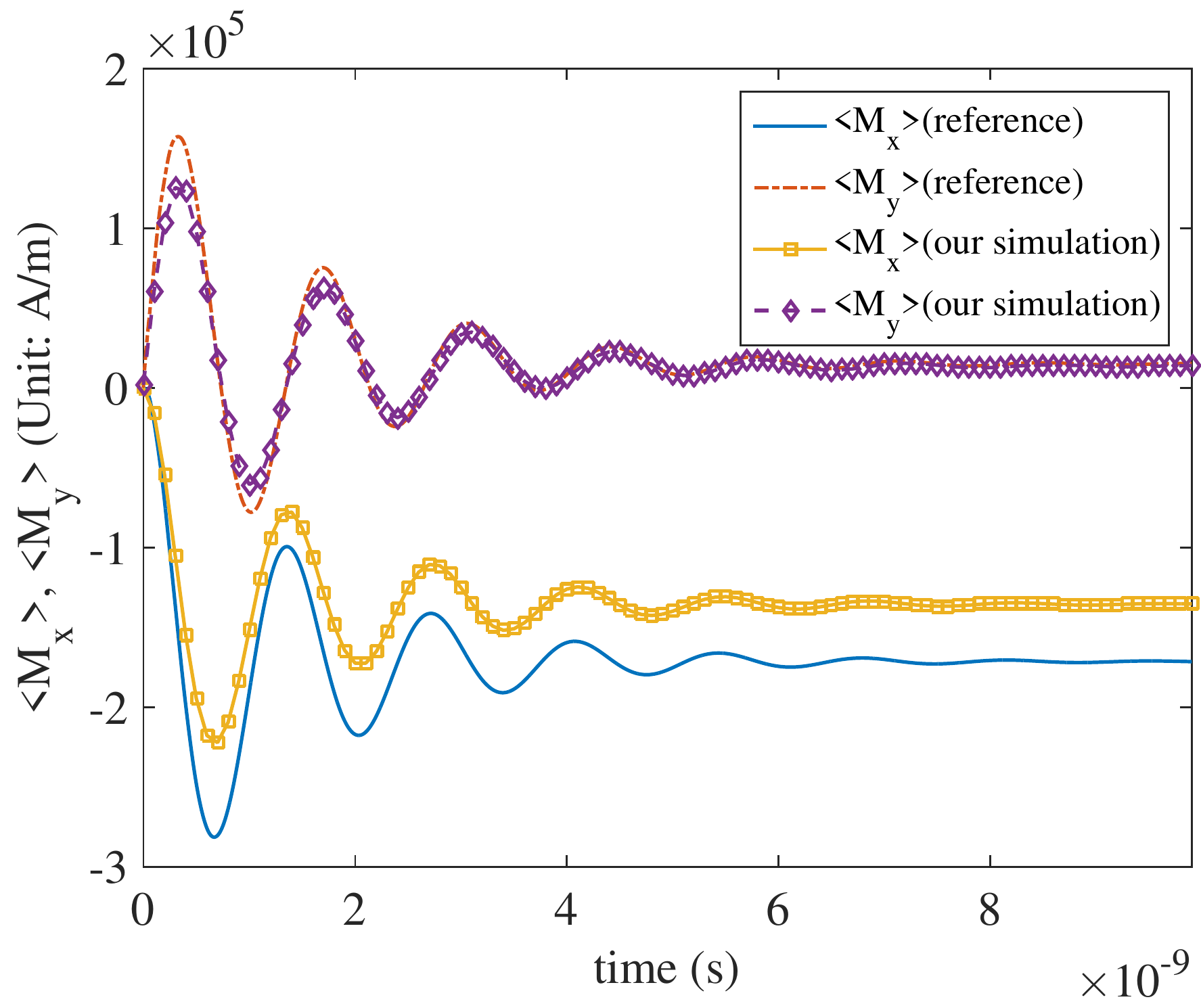}}
		\subfloat[Final state]{\includegraphics[width=2.6in]{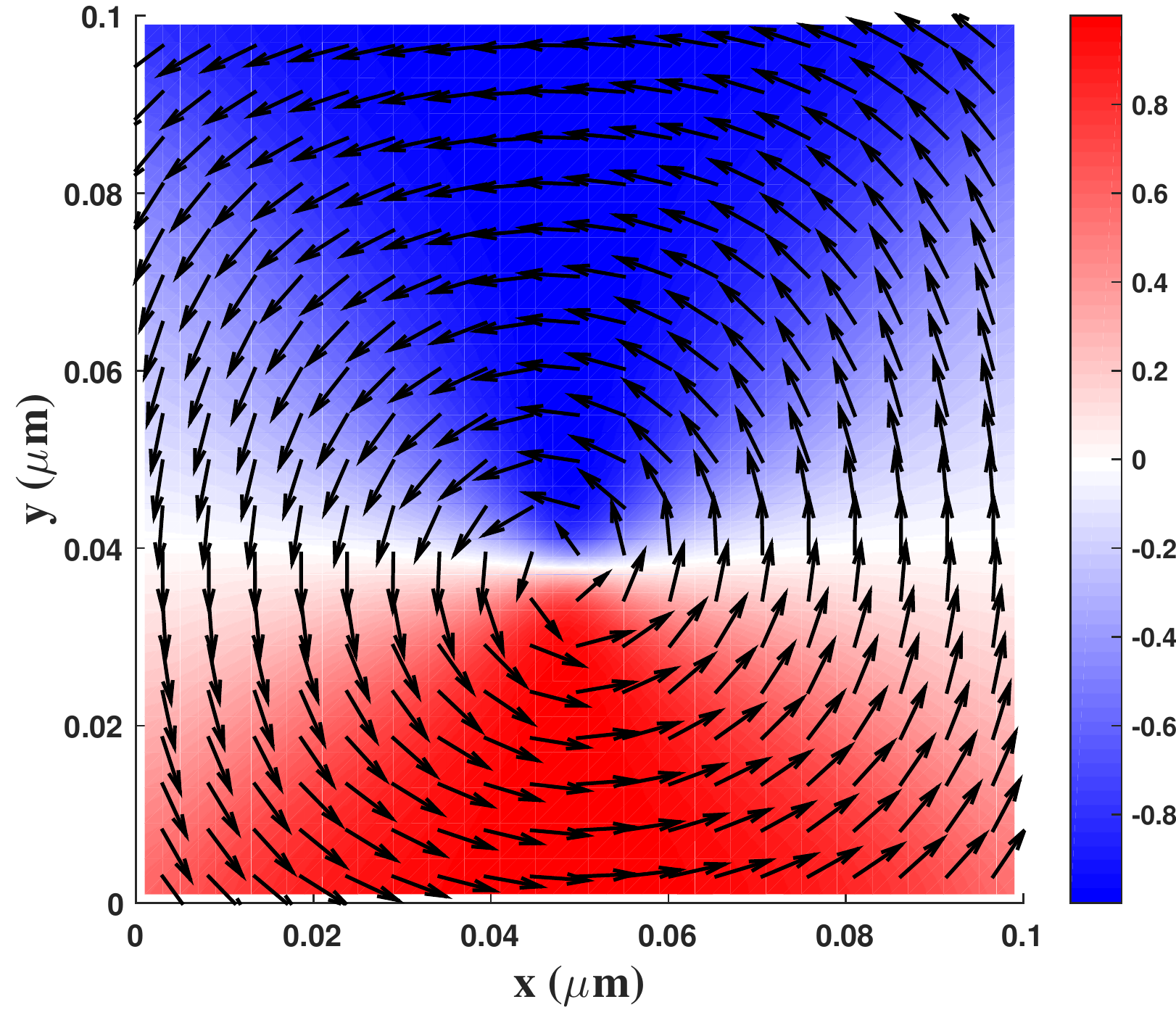}}
		\caption{Magnetization dynamics and the final state in the case of $\xi = 0.05$ and $bJ = 72.17 \;\mathrm{m}/\mathrm{s}$.
			Left: dynamics of the spatially averaged magnetization with the result of D. G. Porter as the reference; Right: the final state colored by the $x$-component of magnetization.}
		\label{fig:xi=0.05}
	\end{figure}
	
	\item[3)] $bJ = 71.64 \mathrm{m}/\mathrm{s}$ and $\xi = 0.1$;
	\begin{figure}[htbp]
		\centering
		\subfloat[Magnetization dynamics]{\includegraphics[width=2.8in]{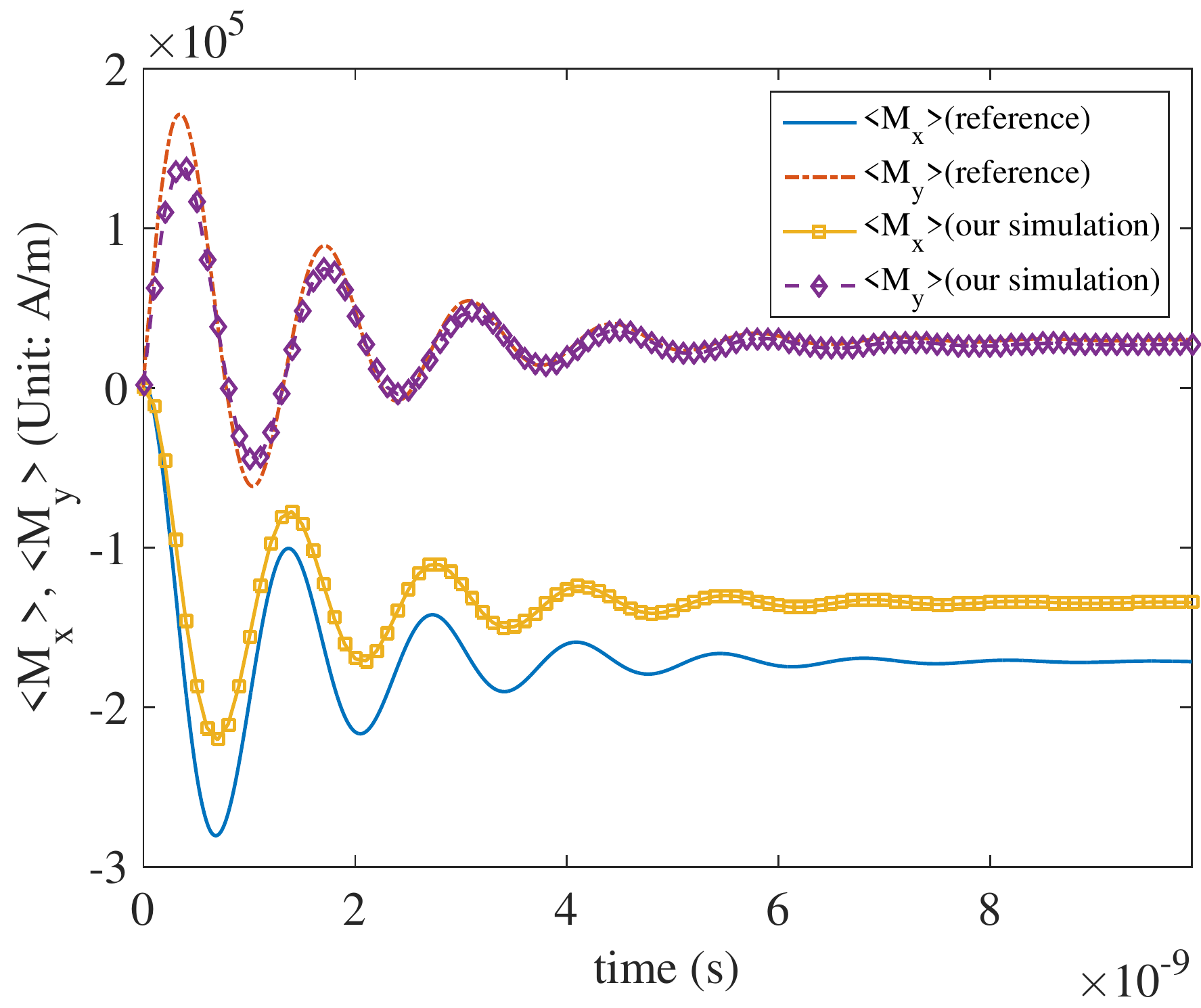}}
		\subfloat[Final state]{\includegraphics[width=2.6in]{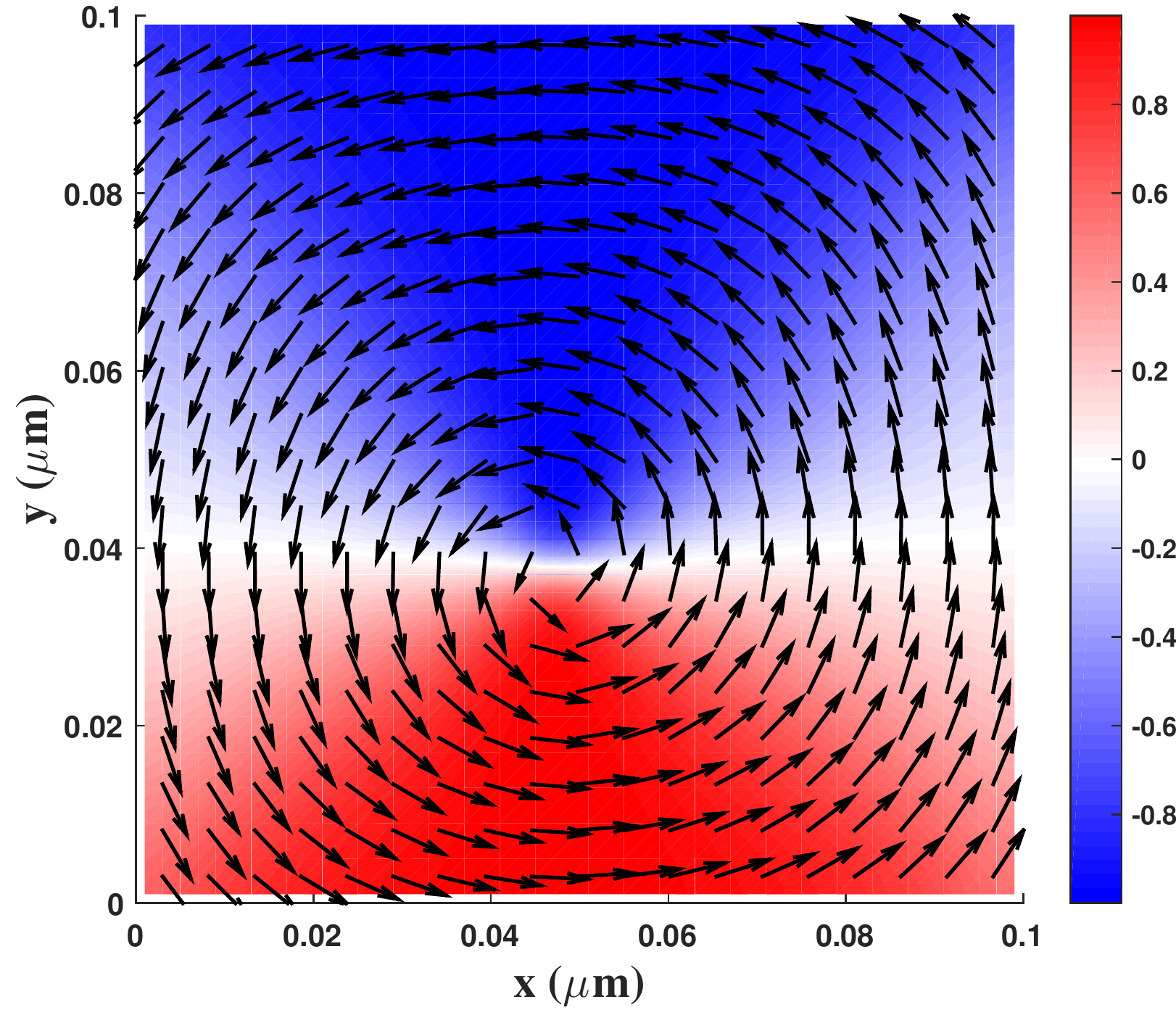}}
		\caption{Magnetization dynamics and the final state in the case of $\xi = 0.1$ and $bJ = 71.64 \;\mathrm{m}/\mathrm{s}$.
			Left: dynamics of the spatially averaged magnetization with the result of D. G. Porter as the reference; Right: the final state colored by the $x$-component of magnetization.}
		\label{fig:xi=0.1}
	\end{figure}
	
	\item[4)] $bJ = 57.88 \mathrm{m}/\mathrm{s}$ and $\xi = 0.5$.
	\begin{figure}[htbp]
		\centering
		\subfloat[Magnetization dynamics]{\includegraphics[width=2.8in]{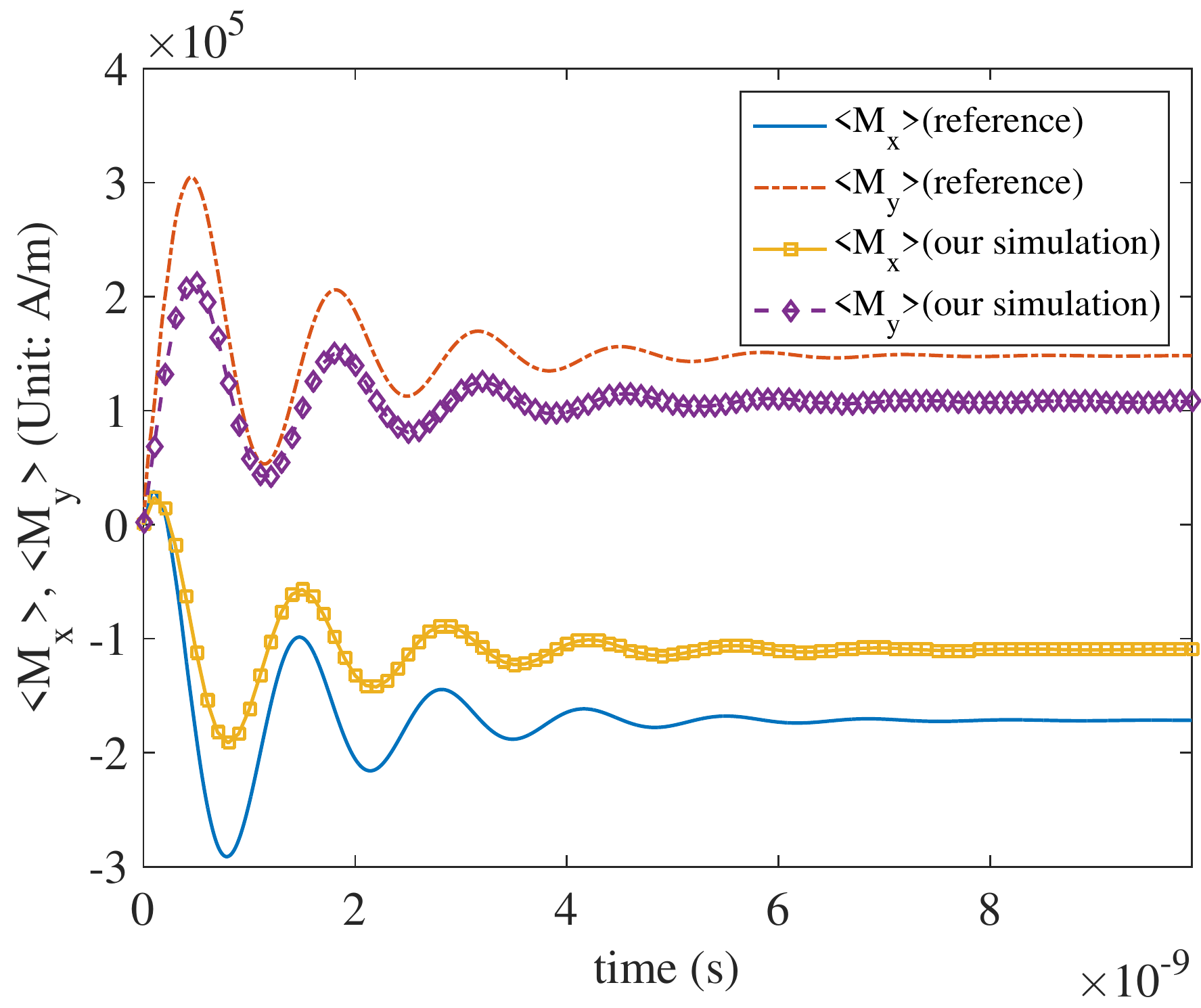}}
		\subfloat[Final state]{\includegraphics[width=2.6in]{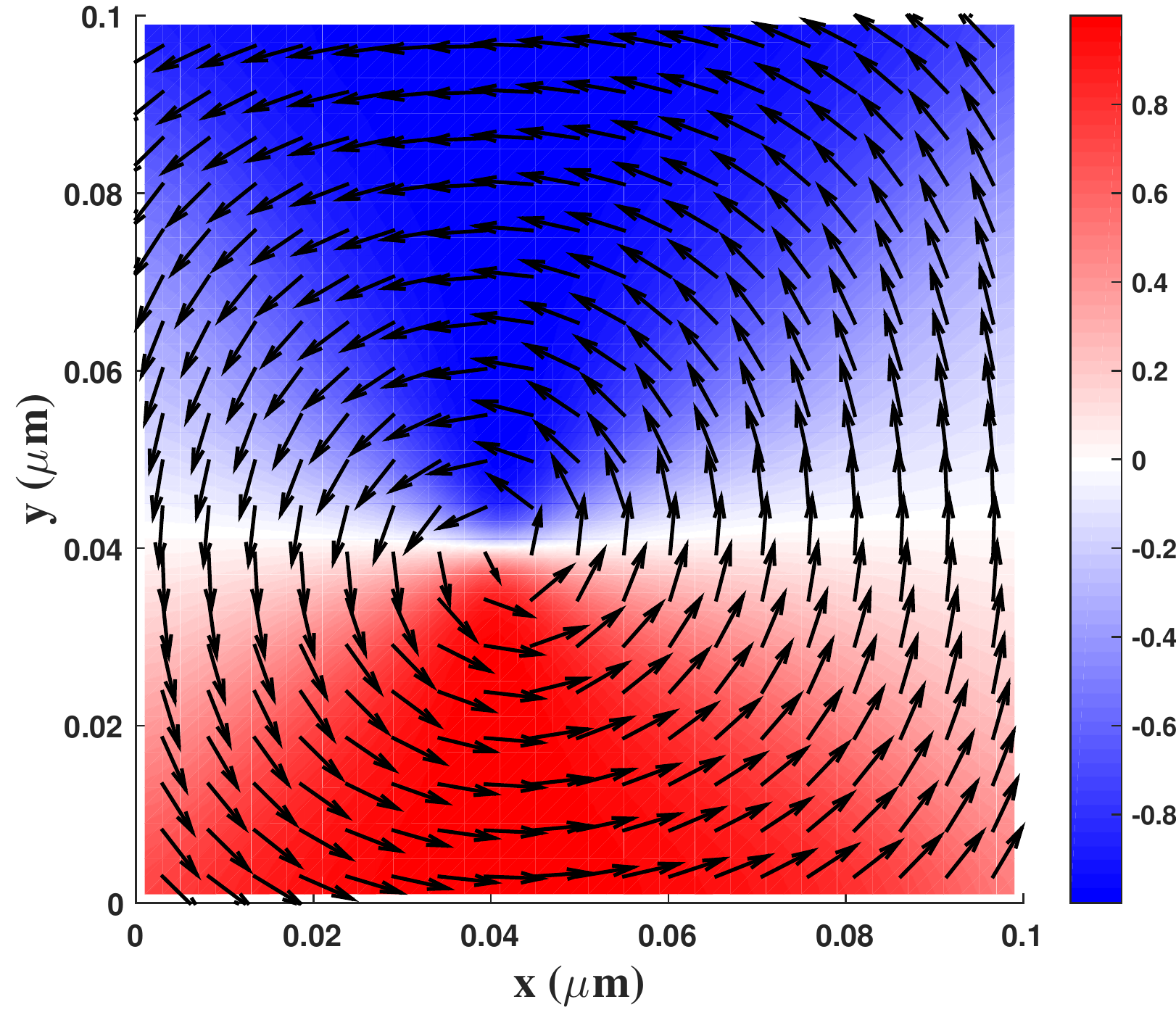}}
		\caption{Magnetization dynamics and the final state in the case of $\xi = 0.5$ and $bJ = 57.88 \;\mathrm{m}/\mathrm{s}$.
			Left: dynamics of the spatially averaged magnetization with the result of D. G. Porter as the reference; Right: the final state colored by the $x$-component of magnetization.}
		\label{fig:xi=0.5}
	\end{figure}
\end{itemize}
In the current work, the first component of the spatially averaged magnetization $<M_x>$ at $10\mathrm{ns}$ is $-1.43\times10^5\;\mathrm{A}/\mathrm{m}$, while the value is  $-1.71\times10^5\;\mathrm{A}/\mathrm{m}$ produced by OOMMF for case 1), 2), and 3). We attribute this discrepancy to different calculations of the stray field. This has been observed in \cite{Najafi2009std} that the spatially averaged magnetization simulated by OOMMF and M$^3$S are different due to different evaluations of the stray field. In our implementation, the stray field is evaluated in 3D without any simplification. For case 4), a vortex state with different location is observed in our simulation. Both $<M_y>$ and $<M_x>$ are different compared to the result of D. G. Porter. However, if the result of G. Finocchio et al. is used for comparison with the parameters $\xi = 0.5$ and $bJ = 72.45 \;\mathrm{m}/\mathrm{s}$, our result has a small difference in $<M_x>$ but no difference in $<M_y>$; as shown in \cref{fig:xi=05andu=72.45}. This state is found to be stable and may be caused by different treatments of STT terms in the LLG equation~\cite{NIST}.

For each case, GSPM-BDF2 takes about $100$ seconds while softwares such as NMAG~\cite{NMAG}, OOMMF~\cite{OOMMF}, and the work in~\cite{Najafi2009std} take more than half an hour on the same personal computer. Below we list the computational costs of GSPM-BDF2 and OOMMF in \cref{tab:ex5}. The proposed method saves $96\%$ computational time over OOMMF.
\begin{table}[htbp]
	\centering
	\caption{Computational costs of the proposed method and OOMMF for Standard Problem \#5~(unit: seconds).}\label{tab:ex5}
	\begin{tabular}{||c|cc|c||}
		\hline
	Standard Problem \#5	& The proposed method & OOMMF & Saving\\
		\hline
		case 1) & 97.58 & 2216.85 & 96\%\\
		case 2) & 97.55 & 2226.45 & 96\%\\
		case 3) & 93.91 & 2229.22 & 96\% \\
		case 4) & 95.59 & 2246.40 & 96\%\\
		\hline
	\end{tabular}
\end{table}

\section{Conclusions}
\label{sec:conclusion}
In this work, we propose a new method (GSPM-BDF2 in short) by combing the advantages of the Gauss-Seidel projection method and the semi-implicit BDF2 scheme. The proposed method solves linear systems of equations with constant coefficients and the spd structure five times and updates the stray field only once. The method is tested to be unconditionally stable with respect to the damping parameter, first-order accurate in time, and second-order accurate in space. In micromagnetic simulations, the proposed method reduces the number of evaluations of the stray field from $3$ to $1$, yielding about $60\%$ reduction of computational time on top of GSPM~\cite{panchi2020GSPM}. For Standard Problem \#4 and \#5 from National Institute of Standards and Technology~\cite{NIST}, GSPM-BDF2 reduces the computational time over OOMMF~\cite{OOMMF} by $82\%$ and $96\%$, respectively. Thus, the proposed method provides a more efficient choice for micromagnetic simulations.

We shall mention that the number of updates for the stray field is minimized in the current work and its evaluation is implemented using FFT which only applies to regular geometries. How to effectively evaluate the stray field over a general geometry is still a difficult task. Realizing this will maximize the applicability of the proposed method for micromagnetic simulations in general.

\section*{Acknowledgments}
P. Li is grateful to Kelong Cheng for helpful discussions during the 18th CSIAM annual meeting and acknowledges the financial support from the Postgraduate Research \& Practice Innovation Program of Jiangsu Province via grant KYCX20\_2711. The work of R. Du was supported in part by NSFC via grant 11501399. The work of J. Chen was supported by NSFC via grant 11971021.

\bibliographystyle{model1-num-names}
\bibliography{refs}

\end{document}